\documentclass[11pt, a4paper]{article}
\usepackage{amsmath}
\usepackage{wasysym}
\usepackage{latexsym}
\usepackage{amsfonts}
\usepackage{mathrsfs}
\usepackage{amssymb}
\usepackage{ifsym}
\usepackage{dsfont}
\usepackage[all]{xy}
\usepackage{float}

\usepackage{siunitx} % Required for alignment

\newtheorem{Definitions1}{Definition}[section]

\newtheorem{Theorems1}{Theorem}[section]

\newtheorem{Coroll1}[Theorems1]{Corollary}
\newtheorem{Lemma1}[Theorems1]{Lemma}

\newtheorem{Examp1}{Example}[section]

\newenvironment{proof}[1][Proof]{\begin{trivlist}
\item[\hskip \labelsep {\bfseries #1}]}{\end{trivlist}}

\newcommand{\qed}{\nobreak \ifvmode \relax \else
      \ifdim\lastskip<1.5em \hskip-\lastskip
      \hskip1.5em plus0em minus0.5em \fi \nobreak
      \vrule height0.75em width0.5em depth0.25em\fi}

\begin{document}
\title{Partially-elementary end extensions of countable models of set theory}
\author{Zachiri McKenzie\\
University of Chester\\
{\tt z.mckenzie@chester.ac.uk}}
\maketitle

\begin{abstract}
Let $\mathsf{KP}$ denote Kripke-Platek Set Theory and let $\mathsf{M}$ be the weak set theory obtained from $\mathsf{ZF}$ by removing the collection scheme, restricting separation to $\Delta_0$-formulae and adding an axiom asserting that every set is contained in a transitive set ($\mathsf{TCo}$). A result due to Kaufmann \cite{kau81} shows that every countable model, $\mathcal{M}$, of $\mathsf{KP}+\Pi_n\textsf{-Collection}$ has a proper $\Sigma_{n+1}$-elementary end extension. We show that for all $n \geq 1$, there exists an $L_\alpha$ that satisfies $\textsf{Separation}$, $\textsf{Powerset}$ and $\Pi_n\textsf{-Collection}$, but that has no $\Sigma_{n+1}$-elementary end extension satisfying either $\Pi_n\textsf{-Collection}$ or $\Pi_{n+3}\textsf{-Foundation}$. Thus showing that there are limits to the amount of the theory of $\mathcal{M}$ that can be transferred to the end extensions that are guaranteed by Kaufmann's Theorem. Using admissible covers and the Barwise Compactness Theorem, we show that if $\mathcal{M}$ is a countable model $\mathsf{KP}+\Pi_n\textsf{-Collection}+\Sigma_{n+1}\textsf{-Foundation}$ and $T$ is a recursive theory that holds in $\mathcal{M}$, then there exists a proper $\Sigma_n$-elementary end extension of $\mathcal{M}$ that satisfies $T$. We use this result to show that the theory $\mathsf{M}+\Pi_n\textsf{-Collection}+\Pi_{n+1}\textsf{-Foundation}$ proves $\Sigma_{n+1}\textsf{-Separation}$.  
\end{abstract}

\section[Introduction]{Introduction}

Keisler and Morley \cite{km68} prove that every countable model of $\mathsf{ZF}$ has a proper elementary end extension. Kaufmann \cite{kau81} refines this result showing that if $n \geq 1$ and $\mathcal{M}$ is a countable structure in the language of set theory that satisfies $\mathsf{KP}+\Pi_n\textsf{-Collection}$, then $\mathcal{M}$ has proper $\Sigma_{n+1}$-elementary end extension. And, conversely, if $n \geq 1$ and $\mathcal{M}$ is a structure in the language of set theory that satisfies $\mathsf{KP}+\mathsf{V=L}$ and has a proper $\Sigma_{n+1}$-elementary end extension, then $\mathcal{M}$ satisfies $\Pi_n\textsf{-Collection}$. A natural question to ask is how much of the theory of $\mathcal{M}$ satisfying $\mathsf{KP}+\Pi_{n}\textsf{-Collection}$ can be made to hold in a proper $\Sigma_{n+1}$-elementary end extension whose existence is guaranteed by Kaufmann's result? In particular, is there a proper $\Sigma_{n+1}$-elementary end extension of $\mathcal{M}$ that also satisfies $\mathsf{KP}+\Pi_{n}\textsf{-Collection}$? Or, if $\mathcal{M}$ is transitive, is there a proper $\Sigma_{n+1}$-elementary end extension of $\mathcal{M}$ that satisfies full induction for all set-theoretic formulae? In section \ref{sec:negativeresults} we show that the answers to the latter two of these questions is ``no". For $n \geq 1$, there is an $L_\alpha$ satisfying $\textsf{Separation}$, $\textsf{Powerset}$ and $\Pi_n\textsf{-Collection}$ that has no proper $\Sigma_{n+1}$-elementary end extension satisfying either $\Pi_n\textsf{-Collection}$ or $\Pi_{n+3}\textsf{-Foundation}$. A key ingredient is a generalisation of a result due to Simpson (see \cite[Remark 2]{kau81}) showing that if $n \geq 1$ and $\mathcal{M}$ is a structure in the language of set theory satisfying $\mathsf{KP}+\mathsf{V=L}$ that has $\Sigma_n$-elementary end extension satisfying enough set theory and with a new ordinal but no least new ordinal, then $\mathcal{M}$ satisfies $\Pi_n\textsf{-Collection}$. Here ``enough set theory" is either $\mathsf{KP}+\Pi_{n-1}\textsf{-Collection}$ or $\mathsf{KP}+\Pi_{n+2}\textsf{-Foundation}$. 

In section \ref{sec:positiveresults}, we use Barwise's admissible cover machinery to build partially-elementary end extensions that satisfy significant fragments of the theory of the model being extended. In particular, we show that if $T$ is a recursively enumerable theory in the language of set theory that extends $\mathsf{KP}+\Pi_n\textsf{-Collection}+\Sigma_{n+1}\textsf{-Foundation}$ and $\mathcal{M}$ is a structure that satisfies $T$, then $\mathcal{M}$ has a proper $\Sigma_n$-elementary end extension that satisfies $T$. That is, by settling for less elementarity we can ensure that there exists an end extension that satisfies any recursively enumerable theory that holds in the model being extended. 

The end-extension result proved in \ref{sec:positiveresults} is used in section \ref{sec:application} to shed light on the relationship between subsystems of $\mathsf{ZF}$ that include the $\textsf{Powerset}$ axiom. We use $\mathsf{M}$ to denote that set theory that is axiomatised by: $\textsf{Extensionality}$, $\textsf{Emptyset}$, $\textsf{Pair}$, $\textsf{Powerset}$, $\textsf{TCo}$, $\textsf{Infinity}$, $\Delta_0\textsf{-Separation}$ and $\textsf{Set-Foundation}$. We show that for all $n \geq 1$, $\mathsf{M}+\Pi_n\textsf{-Collection}+\Pi_{n+1}\textsf{-Foundation}$ proves $\Sigma_{n+1}\textsf{-Separation}$. In particular, for all $n \geq 1$, the theories $\mathsf{M}+\Pi_n\textsf{-Collection}$ and $\mathsf{M}+\textsf{Strong } \Pi_n\textsf{-Collection}$ have the same well-founded models, settling a question about heights of minimum models of subsystems of $\mathsf{ZF}$ including $\textsf{Powerset}$ left open in \cite{gos80}.

\section[Background]{Background} \label{sec:background}

Let $\mathcal{L}$ be the language of set theory-- the language whose only non-logical symbol is the binary relation $\in$. Let $\mathcal{L}^\prime$ be a language that contains $\mathcal{L}$ and let $\Gamma$ be a collection of $\mathcal{L}^\prime$-formulae.
\begin{itemize}
\item $\Gamma\textsf{-Separation}$ is the scheme that consists of the sentences 
$$\forall \vec{z} \forall w \exists y \forall x (x \in y \iff (x \in w \land \phi(x, \vec{z})),$$
for all formulae $\phi(x, \vec{z})$ in $\Gamma$. $\textsf{Separation}$ is the scheme that consists of these sentences for every formula $\phi(x, \vec{z})$ in $\mathcal{L}$.
\item $\Gamma\textsf{-Collection}$ is the scheme that consists of the sentences
\[
\forall \vec{z} \forall w ((\forall x \in w) \exists y \phi(x, y, \vec{z}) \Rightarrow \exists c (\forall x \in w) (\exists y \in c)\phi(x, y, \vec{z})),
\]
for all formulae $\phi(x, y, \vec{z})$ in $\Gamma$. $\textsf{Collection}$ is the scheme that consists of these sentences for every formula $\phi(x, y, \vec{z})$ in $\mathcal{L}$.
\item $\textsf{Strong }\Gamma\textsf{-Collection}$ is the scheme that consists of the sentences
\[
\forall \vec{z} \forall w \exists c (\forall x \in w)( \exists y \phi(x, y, \vec{z}) \Rightarrow (\exists y \in c) \phi(x, y, \vec{z})),
\]
for all formulae $\phi(x, y, \vec{z})$ in $\Gamma$. $\textsf{Strong Collection}$ is the scheme that consists of these sentences for every formula $\phi(x, y, \vec{z})$ in $\mathcal{L}$.
\item $\Gamma\textsf{-Foundation}$ is the scheme that consists of the sentences 
\[
\forall \vec{z} (\exists x \phi(x, \vec{z}) \Rightarrow \exists y (\phi(y, \vec{z}) \land (\forall w \in y) \neg \phi(w, \vec{z}))),
\]
for all formulae $\phi(x, \vec{z})$ in $\Gamma$. If $\Gamma= \{x \in z\}$, then the resulting axiom is referred to as $\textsf{Set-Foundation}$. $\textsf{Foundation}$ is the scheme that consists of these sentences for every formula $\phi(x, \vec{z})$ in $\mathcal{L}$.  
\end{itemize}
In addition to the L\'{e}vy classes of $\mathcal{L}$-formulae, $\Delta_0$, $\Sigma_1$, $\Pi_1$, \ldots, we will also make reference to the class $\Delta_0^\mathcal{P}$, introduced in \cite{tak72}, that consists of $\mathcal{L}$-formulae whose quantifiers are bounded either by the membership relation ($\in$) or the subset relation ($\subseteq$), and the classes $\Sigma_1^\mathcal{P}$, $\Pi_1^\mathcal{P}$, $\Sigma_2^\mathcal{P}$, \ldots that are defined from $\Delta_0^\mathcal{P}$ in the same way that the classes $\Sigma_1$, $\Pi_1$, $\Sigma_2$, \ldots are defined from $\Delta_0$. Let $T$ be a theory in a language, $\mathcal{L}^\prime$, that includes $\mathcal{L}$. Let $\Gamma$ be a class of $\mathcal{L}^\prime$-formulae. A formula is $\Gamma$ in $T$ or $\Gamma^T$ if it is provably equivalent in $T$ to a formula in $\Gamma$. A formula is $\Delta_n$ in $T$ or $\Delta_n^T$ if it is both $\Sigma_n^T$ and $\Pi_n^T$. 
\begin{itemize}
\item $\Delta_n\textsf{-Separation}$ is the scheme that consists of the sentences 
$$\forall \vec{z}(\forall v (\phi(v, \vec{z}) \iff \psi(v, \vec{z})) \Rightarrow \forall w \exists y \forall x(x \in y \iff (x \in w \land \phi(x, \vec{z}))))$$
for all $\Sigma_n$-formulae $\phi(x, \vec{z})$ and $\Pi_n$-formulae $\psi(x, \vec{z})$.
\item $\Delta_n\textsf{-Foundation}$ is the scheme that consists of the sentences 
$$\forall \vec{z}(\forall v(\phi(x, \vec{z}) \iff \psi(x, \vec{z})) \Rightarrow (\exists x \phi(x, \vec{z}) \Rightarrow \exists y (\phi(y, \vec{z}) \land (\forall w \in y) \neg \phi(w, \vec{z}))))$$
for all $\Sigma_n$-formulae $\phi(x, \vec{z})$ and $\Pi_n$-formulae $\psi(x, \vec{z})$. 
\end{itemize} 
We use $\mathsf{S}_1$ to denote the $\mathcal{L}$-theory with axioms: \textsf{Extensionality}, \textsf{Emptyset}, \textsf{Pair}, \textsf{Union}, \textsf{Set Difference} and \textsf{Powerset}. Following \cite{mat01}, we take {\bf Kripke-Platek Set Theory} ($\mathsf{KP}$) to be the theory obtained from $\mathsf{S}_1$ by removing \textsf{Powerset} and adding $\Delta_0\textsf{-Separation}$, $\Delta_0\textsf{-Collection}$ and $\Pi_1\textsf{-Foundation}$. Note that this differs from \cite{bar75, fri73}, which defines Kripke-Platek Set Theory to include \textsf{Foundation}. The theory $\mathsf{KPI}$ is obtained from $\mathsf{KP}$ by adding the axiom \textsf{Infinity}, which states that a superset of the von Neumann ordinal $\omega$ exists. We use $\mathsf{M}^-$ to denote the theory that is obtained from $\mathsf{KPI}$ by replacing $\Pi_1\textsf{-Foundation}$ with $\textsf{Set-Foundation}$ and removing $\Delta_0\textsf{-Collection}$, and adding an axiom $\mathsf{TCo}$ asserting that every set is contained in a transitive set. The theory $\mathsf{M}$ is obtained from $\mathsf{M}^-$ by adding $\textsf{Powerset}$. The theory $\mathsf{MOST}$ is obtained form $\mathsf{M}$ by adding $\textsf{Strong } \Delta_0\textsf{-Collection}$ and the Axiom of Choice ($\mathsf{AC}$). {\bf Zermelo Set Theory} ($\mathsf{Z}$) is obtained for $\mathsf{M}$ by removing $\mathsf{TCo}$ and adding $\textsf{Separation}$. The theory $\mathsf{KP}^\mathcal{P}$ is obtained from $\mathsf{M}$ by adding $\Delta_0^\mathcal{P}\textsf{-Collection}$ and $\Pi_1^\mathcal{P}\textsf{-Foundation}$. 

The theory $\mathsf{KP}$ proves $\mathsf{TCo}$ (see, for example, \cite[I.6.1]{bar75}). Both $\mathsf{KP}$ and $\mathsf{M}$ prove that every set $x$ is contained in a least transitive set that is called the {\bf transitive closure} of $x$, and denoted $\mathsf{TC}(x)$. The following are some important consequences of fragments of the collection scheme over the theory $\mathsf{M}^-$:
\begin{itemize}
\item The proof of \cite[I.4.4]{bar75} generalises to show that, in the theory $\mathsf{M}^-$, $\Pi_n\textsf{-Collection}$ implies $\Sigma_{n+1}\textsf{-Collection}$.
\item \cite[Lemma 4.13]{flw16} shows that, over $\mathsf{M}^-$, $\Pi_n\textsf{-Collection}$ implies $\Delta_{n+1}\textsf{-Separation}$.
\item It is noted in \cite[Proposition 2.4]{flw16} that if $T$ is $\mathsf{M}^-+\Pi_n\textsf{-Collection}$, then the classes $\Sigma_{n+1}^T$ and $\Pi_{n+1}^T$ are closed under bounded quantification.
\item \cite[Lemma 2.4]{mck19}, for example, shows that, over $\mathsf{M}^-$, $\textsf{Strong }\Pi_n\textsf{-Collection}$ is equivalent to $\Pi_n\textsf{-Collection}+\Sigma_{n+1}\textsf{-Separation}$. 
\end{itemize} 

Let $\mathcal{L}^\prime$ be a language that contains $\mathcal{L}$. Let $\mathcal{M}= \langle M, \in^\mathcal{M}, \ldots \rangle$ be an $\mathcal{L}^\prime$-structure. If $a \in M$, then we will use $a^*$ to denote the set $\{x \in M \mid \mathcal{M} \models (x \in a)\}$, as long as $\mathcal{M}$ is clear from the context. Let $\Gamma$ be a collection of $\mathcal{L}^\prime$-formulae. We say $X \subseteq M$ is $\Gamma$ over $\mathcal{M}$ if there is a formula $\phi(x, \vec{z})$ in $\Gamma$ and $\vec{a} \in M$ such that $X= \{x \in M \mid \mathcal{M} \models \phi(x, \vec{a})\}$. In the special case that $\Gamma$ is all $\mathcal{L}^\prime$-formulae, we say that $X$ is a {\bf definable subclass} of $\mathcal{M}$. A set $X \subseteq M$ is $\Delta_n$ over $\mathcal{M}$ if it is both $\Sigma_n$ over $\mathcal{M}$ and $\Pi_n$ over $\mathcal{M}$. 

A structure $\mathcal{N}= \langle N, \in^\mathcal{N} \rangle$ is an {\bf end extension} of $\mathcal{M}= \langle M, \in^\mathcal{M} \rangle$, written $\mathcal{M} \subseteq_e \mathcal{N}$, if $\mathcal{M}$ is a substructure of $\mathcal{N}$ and for all $x \in M$ and for all $y \in N$, if $\mathcal{N} \models (y \in x)$, then $y \in M$. An end extension $\mathcal{N}$ of $\mathcal{M}$ is {\bf proper} if $M \neq N$. If $\mathcal{N}= \langle N, \in^\mathcal{N} \rangle$ is an end extension of $\mathcal{M}= \langle M, \in^\mathcal{M} \rangle$ and for all $x \in M$ and for all $y \in N$, if $\mathcal{N} \models (y \subseteq x)$, then $y \in M$, then we say that $\mathcal{N}$ is a {\bf powerset-preserving end extension} of $\mathcal{M}$ and write $\mathcal{M} \subseteq_e^\mathcal{P} \mathcal{N}$. We say that $\mathcal{N}$ is a {\bf $\Sigma_n$-elementary end extension} of $\mathcal{M}$, and write $\mathcal{M} \prec_{e, n} \mathcal{N}$, if $\mathcal{M} \subseteq_e \mathcal{N}$ and $\Sigma_n$ properties are preserved between $\mathcal{M}$ and $\mathcal{N}$.     

We use $\mathsf{Ord}$ to denote the class of ordinals. As shown in \cite[Chapter II]{bar75}, the theory $\mathsf{KP}$ is capable of defining G\"{o}del's constructible universe ($L$). For all sets $X$,
\[
\mathsf{Def}(X)= \{Y \subseteq X \mid Y \textrm{ is a definable subclass of } \langle X, \in \rangle \},
\]
which can be seen to be a set in the theory $\mathsf{KP}$ using a formula for satisfaction in set structures such as the one described in \cite[Section III.1]{bar75}. The levels of $L$ are then defined by the recursion:
$$L_0= \emptyset \textrm{ and } L_\alpha= \bigcup_{\beta < \alpha} L_\beta \textrm{ if } \alpha \textrm{ is a limit ordinal,}$$
$$L_{\alpha+1}= L_\alpha \cup \mathsf{Def}(L_\alpha), \textrm{ and}$$
$$L= \bigcup_{\alpha \in \mathsf{Ord}} L_\alpha.$$
The function $\alpha \mapsto L_\alpha$ is total in $\mathsf{KP}$ and $\Delta_1^{\mathsf{KP}}$. The axiom $\mathsf{V=L}$ asserts that every set is the member of some $L_\alpha$. A transitive set $M$ such that $\langle M, \in \rangle$ satisfies $\mathsf{KP}$ is said to be an {\bf admissible set}. An ordinal $\alpha$ is said to be an {\bf admissible ordinal} if $L_\alpha$ is an admissible set.

The theory $\mathsf{KP}^\mathcal{P}$ proves that the function $\alpha \mapsto V_\alpha$ is total and $\Delta_1^\mathcal{P}$. Mathias \cite[Proposition Scheme 6.12]{mat01} refines the relationships between the classes $\Delta_0^\mathcal{P}$, $\Sigma_1^\mathcal{P}$, $\Pi_1^{\mathcal{P}}$, \ldots, and the L\'{e}vy classes by showing that $\Sigma_1 \subseteq (\Delta_1^\mathcal{P})^{\textsf{MOST}}$ and $\Delta_0^\mathcal{P} \subseteq \Delta_2^{\mathsf{S}_1}$. Therefore, the function $\alpha \mapsto V_\alpha$ is $\Delta_2^{\mathsf{KP}^\mathcal{P}}$. It also follows from this analysis that $\mathsf{KP}^\mathcal{P}$ is a subtheory of $\mathsf{M}+\Pi_1\textsf{-Collection}+\Pi_2\textsf{-Foundation}$. 

Let $T$ be an $\mathcal{L}$-theory. A transitive set $M$ is said to be the {\bf minimum model} of $T$ if $\langle M, \in \rangle \models T$ and for all transitive sets $N$ with $\langle N, \in \rangle \models T$, $M \subseteq N$. For example, $L_{\omega_1^\mathrm{CK}}$ is the minimum model of $\mathsf{KPI}$. For an $\mathcal{L}$-theory $T$ to have a minimum model it is sufficient that the conjunction of the following conditions hold:
\begin{itemize}
\item[(I)] There exists a transitive set $M$ such that $\langle M, \in \rangle \models T$;
\item[(II)] for all transitive $M$ with $\langle M, \in \rangle \models T$, $\langle L^M, \in \rangle \models T$. 
\end{itemize}
Gostanian \cite[\S 1]{gos80} shows that all sufficiently strong subsystems of $\mathsf{ZF}$ and $\mathsf{ZF}^-$ obtained by restricting the separation and collection schemes to formulae in the L\'{e}vy classes have minimum models. In particular:

\begin{Theorems1}
(Gostanian \cite{gos80}) Let $n, m \in \omega$.
\begin{itemize}
\item[(I)] The theory $\mathsf{KPI}+\Pi_m\textsf{-Separation}+\Pi_n\textsf{-Collection}$ has a minimum model. Moreover, the minimum model of this theory satisfies $\mathsf{V=L}$.
\item[(II)] If $n \geq 1$ or $m \geq 1$, then the theory $\mathsf{KPI}+\textsf{Powerset}+\Pi_m\textsf{-Separation}+\Pi_n\textsf{-Collection}$ has a minimum model. Moreover, the minimum model of this theory satisfies $\mathsf{V=L}$.   
\end{itemize}
\Square
\end{Theorems1}

Gostanian's analysis also yields:

\begin{Theorems1}
Let $n \in \omega$. The theory $\mathsf{Z}+\Pi_n\textsf{-Collection}$ has a minimum model. Moreover, the minimum model of this theory satisfies $\mathsf{V=L}$. \Square 
\end{Theorems1}   

The fact that $\mathsf{KP}$ is able to define satisfaction in set structures also facilitates the definition of formulae expressing satisfaction, in the universe, for formulae in any given level of the L\'{e}vy hierarchy.

\begin{Definitions1} \label{Df:Delta0Satisfaction}
The formula $\mathsf{Sat}_{\Delta_0}(q, x)$ is defined as
$$\begin{array}{c}
(q \in \omega) \land (q= \ulcorner \phi(v_1, \ldots, v_m) \urcorner \textrm{ where } \phi \textrm{ is } \Delta_0) \land (x= \langle x_1, \ldots, x_m \rangle) \land\\
\exists N \left( \bigcup N \subseteq N \land (x_1, \ldots, x_m \in N) \land (\langle N, \in \rangle \models \phi[x_1, \ldots, x_m]) \right)
\end{array}.$$
\end{Definitions1}

We can now inductively define formulae $\mathsf{Sat}_{\Sigma_n}(q, x)$ and $\mathsf{Sat}_{\Pi_n}(q, x)$ that express satisfaction for formulae in the classes $\Sigma_n$ and $\Pi_n$.

\begin{Definitions1}
The formulae $\mathsf{Sat}_{\Sigma_n}(q, x)$ and $\mathsf{Sat}_{\Pi_n}(q, x)$ are defined recursively for $n>0$. $\mathsf{Sat}_{\Sigma_{n+1}}(q, x)$ is defined as the  formula
$$\exists \vec{y} \exists k \exists b \left( \begin{array}{c}
(q= \ulcorner\exists \vec{u} \phi(\vec{u}, v_1, \ldots, v_l)\urcorner \textrm{ where } \phi \textrm{ is } \Pi_n)\land (x= \langle x_1, \ldots, x_l \rangle)\\
\land (b= \langle \vec{y}, x_1, \ldots, x_l \rangle) \land (k= \ulcorner \phi(\vec{u}, v_1, \ldots, v_l) \urcorner) \land \mathsf{Sat}_{\Pi_n}(k, b)
\end{array}\right);$$
and  $\mathsf{Sat}_{\Pi_{n+1}}(q, x)$ is defined as the formula
$$\forall \vec{y} \forall k \forall b \left( \begin{array}{c}
(q= \ulcorner\forall \vec{u} \phi(\vec{u}, v_1, \ldots, v_l) \urcorner \textrm{ where } \phi \textrm{ is } \Sigma_n)\land (x= \langle x_1, \ldots, x_l \rangle)\\
\land ((b= \langle \vec{y}, x_1, \ldots, x_l \rangle) \land (k= \ulcorner\phi(\vec{u}, v_1, \ldots, v_l)\urcorner) \Rightarrow \mathsf{Sat}_{\Sigma_n}(k, b))
\end{array}\right).$$
\end{Definitions1}

\begin{Theorems1} \label{Complexityofpartialsat} 
Suppose $n \in \omega$ and $m=\max \{ 1, n \}$. The formula $\mathsf{Sat}_{\Sigma_n}(q, x)$ (respectively $\mathsf{Sat}_{\Pi_n}(q, x)$) is $\Sigma_m^{\mathsf{KP}}$ ($\Pi_m^{\mathsf{KP}}$, respectively). Moreover, $\mathsf{Sat}_{\Sigma_n}(q, x)$ (respectively $\mathsf{Sat}_{\Pi_n}(q, x)$) expresses satisfaction for $\Sigma_n$-formulae ($\Pi_n$-formulae, respectively) in the theory $\mathsf{KP}$, i.e., if $\mathcal{M} \models\mathsf{KP}$, $\phi(v_1,\ldots,v_k)$ is a $\Sigma_n$-formula, and $x_1,\ldots,x_k$ are in $M$, then for $q = \ulcorner   \phi( v_1, \ldots, v_k) \urcorner$, $\mathcal{M}$ satisfies the universal generalisation of the following formula:
$$x= \langle x_1, \ldots,x_k \rangle \Rightarrow \left( \phi(x_1,\ldots,x_k) \iff \mathrm{Sat}_{\Sigma_{n}}(q, x) \right).$$
\Square
\end{Theorems1}

Kaufmann \cite{kau81} identifies necessary and sufficient conditions for models of $\mathsf{KP}$ to have proper $\Sigma_n$-elementary end extensions.

\begin{Theorems1} \label{Th:KaufmannEEResult} 
(Kaufmann \cite[Theorem 1]{kau81}) Let $n \geq 1$. Let $\mathcal{M}= \langle M, \in^\mathcal{M} \rangle$ be a model of $\mathsf{KP}$. Consider
\begin{itemize}
\item[(I)] there exists $\mathcal{N}= \langle N, \in^\mathcal{N} \rangle$ such that $\mathcal{M} \prec_{e, n+1} \mathcal{N}$ and $M \neq N$;
\item[(II)] $\mathcal{M} \models \Pi_{n}\textsf{-Collection}$.  
\end{itemize}
If $\mathcal{M} \models \mathsf{V=L}$, then $(I) \Rightarrow (II)$. If $M$ is countable, then $(II) \Rightarrow (I)$. \Square   
\end{Theorems1}

It should be noted that Kaufmann proves that (II) implies (I) in the above under the weaker assumption that $\mathcal{M}$ is a resolvable model of $\mathsf{M}^-$. A model $\mathcal{M}=\langle M, \in^\mathcal{M} \rangle$ of $\mathsf{M}^-$ is {\bf resolvable} if there is a function $F$ that is $\Delta_1$ over $\mathcal{M}$ such that for all $x \in M$, there exists $\alpha \in \mathsf{Ord}^\mathcal{M}$ such that $x \in F(\alpha)$. The function $\alpha \mapsto L_\alpha$ witnesses the fact that every model of $\mathsf{KP}+\mathsf{V=L}$ is resolvable.

\section[Limitations of Kaufmann's Theorem]{Limitations of Kaufmann's Theorem} \label{sec:negativeresults}

In this section we show that there are limitations on the amount of the theory of the base model that can be transferred to the partially-elementary end extension guaranteed by Theorem \ref{Th:KaufmannEEResult}. We utilise a generalisation of a result, due to Simpson and that that is mentioned in \cite[Remark 2]{kau81}, showing that if a $\mathcal{M}$ satisfies $\mathsf{KP}+\mathsf{V=L}$ and has a $\Sigma_n$-elementary end extension that satisfies enough set theory and contains no least new ordinal, then $\mathcal{M}$ must satisfy $\Pi_{n}\textsf{-Collection}$. The proof of this generalisation, Theorem \ref{Th:CollectionFromPartiallyElementaryEE1}, is based on Enayat's proof of a refinement of Simpson's result (personal communication) that corresponds to the specific case where $n=1$ and $\mathcal{M}$ is transitive.

\begin{Theorems1} \label{Th:CollectionFromPartiallyElementaryEE1}
Let $n \geq 1$. Let $\mathcal{M}= \langle M, \in^\mathcal{M} \rangle$ be a model of $\mathsf{KP}+\mathsf{V=L}$. Suppose $\mathcal{N}= \langle N, \in^\mathcal{N} \rangle$ is such that $\mathcal{M} \prec_{e, n} \mathcal{N}$, $\mathcal{N} \models \mathsf{KP}$ and $\mathsf{Ord}^\mathcal{N} \backslash \mathsf{Ord}^\mathcal{M}$ is nonempty and has no least element. If $\mathcal{N}\models \Pi_{n-1}\textsf{-Collection}$ or $\mathcal{N} \models \Pi_{n+2}\textsf{-Foundation}$, then $\mathcal{M} \models \Pi_n\textsf{-Collection}$.
\end{Theorems1}

\begin{proof}
Assume that $\mathcal{N}= \langle N, \in^\mathcal{N} \rangle$ is such that 
\begin{itemize}
\item[(I)] $\mathcal{M} \prec_{e, n} \mathcal{N}$;
\item[(II)] $\mathcal{N} \models \mathsf{KP}$;
\item[(III)] $\mathsf{Ord}^\mathcal{N} \backslash \mathsf{Ord}^\mathcal{M}$ is nonempty and has no least element. 
\end{itemize}
Note that, since $\mathcal{M} \prec_{e, 1} \mathcal{N}$ and $\mathcal{M} \models \mathsf{V=L}$, for all $\beta \in \mathsf{Ord}^\mathcal{N} \backslash \mathsf{Ord}^\mathcal{M}$, $M \subseteq (L_\beta^{\mathcal{N}})^*$. We need to show that if either $\Pi_{n-1}\textsf{-Collection}$ or $\Pi_{n+2}\textsf{-Foundation}$ hold in $\mathcal{N}$, then $\mathcal{M} \models \Pi_n\textsf{-Collection}$. Let $\phi(x, y, \vec{z})$ be a $\Pi_n$-formula. Let $\vec{a}, b \in M$ be such that 
\[
\mathcal{M} \models (\forall x \in b) \exists y \phi(x, y, \vec{a}).
\]
So, for all $x \in b^*$, there exists $y \in M$ such that
$$\mathcal{M} \models \phi(x, y, \vec{a}).$$
Therefore, since $\mathcal{M} \prec_{e, n} \mathcal{N}$, for all $x \in b^*$, there exists $y \in M$ such that
$$\mathcal{N} \models \phi(x, y, \vec{a}).$$
Now, $\phi(x, y, \vec{z})$ can be written as $\forall w \psi(w, x, y, \vec{z})$ where $\psi(w, x, y, \vec{z})$ is $\Sigma_{n-1}$. Let $\xi \in \mathsf{Ord}^\mathcal{N} \backslash \mathsf{Ord}^\mathcal{M}$. So, for all $\beta \in \mathsf{Ord}^\mathcal{N} \backslash \mathsf{Ord}^\mathcal{M}$ and for all $x \in b^*$, there exists $y \in (L_\beta^\mathcal{N})^*$ such that
$$\mathcal{N} \models (\forall w \in L_\xi)\psi(w, x, y, \vec{a}).$$
Therefore, for all $\beta \in \mathsf{Ord}^\mathcal{N} \backslash \mathsf{Ord}^\mathcal{M}$,
\begin{equation}\label{eq:Equation1}
\mathcal{N} \models (\forall x \in b)(\exists y \in L_\beta)(\forall w \in L_\xi) \psi(w, x, y, \vec{a})
\end{equation}
Now, define $\theta(\beta, \xi, b, \vec{a})$ to be the formula
$$(\forall x \in b)(\exists y \in L_\beta)(\forall w \in L_\xi) \psi(w, x, y, \vec{a}).$$
If $\Pi_{n-1}\textsf{-Collection}$ holds in $\mathcal{N}$, then $\theta(\beta, \xi, b, \vec{a})$ is equivalent to a $\Sigma_{n-1}$-formula. Without $\Pi_{n-1}\textsf{-Collection}$, $\theta(\beta, \xi, b, \vec{a})$ can be written as a $\Pi_{n+2}$-formula. Therefore, $\Pi_{n-1}\textsf{-Collection}$ or $\Pi_{n+2}\textsf{-Foundation}$ in $\mathcal{N}$ will ensure that there is a least $\beta_0 \in \mathrm{Ord}^\mathcal{N}$ such that $\mathcal{N} \models \theta(\beta_0, \xi, b, \vec{a})$. Moreover, by (\ref{eq:Equation1}), $\beta_0 \in M$. Therefore,
$$\mathcal{N} \models (\forall x \in b) (\exists y \in L_{\beta_0})(\forall w \in L_\xi) \psi(w, x, y, \vec{a}).$$
So, for all $x \in b^*$, there exists $y \in (L_{\beta_0}^\mathcal{M})^*$, for all $w \in (L_\xi^\mathcal{N})^*$, 
$$\mathcal{N} \models \psi(w, x, y, \vec{a}).$$
Which, since $\mathcal{M} \prec_{e, n} \mathcal{N}$, implies that for all $x \in b^*$, there exists $y \in (L_{\beta_0}^\mathcal{M})^*$, for all $w \in M$, 
$$\mathcal{M} \models \psi(w, x, y, \vec{a}).$$
Therefore, $\mathcal{M} \models (\forall x \in b)(\exists y \in L_{\beta_0}) \phi(x, y, \vec{a})$. This shows that $\Pi_n\textsf{-Collection}$ holds in $\mathcal{M}$.    
\Square
\end{proof}  

Enayat uses a specific case of Theorem \ref{Th:CollectionFromPartiallyElementaryEE1} to show that the \mbox{$\langle L_{\omega_1^{\mathrm{CK}}}, \in \rangle$} has no proper $\Sigma_1$-elementary end extension that satisfies $\mathsf{KP}$. We now turn to generalising this result to show that for all $n \geq 1$, the minimum model of $\mathsf{Z}+\Pi_n\textsf{-Collection}$ has no proper $\Sigma_{n+1}$-elementary end extension that satisfies either $\mathsf{KP}+\Pi_{n+3}\textsf{-Foundation}$ or $\mathsf{KP}+\Pi_n\textsf{-Collection}$. However, by Theorem \ref{Th:KaufmannEEResult}, for all $n \geq 1$, the minimum model of $\mathsf{Z}+\Pi_n\textsf{-Collection}$ does have a proper $\Sigma_{n+1}$-elementary end extension.

The following result follows from \cite[Theorem 4.4]{mck19}:

\begin{Theorems1}
Let $n \geq 1$. The theory $\mathsf{M}+\Pi_{n+1}\textsf{-Collection}+\Pi_{n+2}\textsf{-Foundation}$ proves that there exists a transitive model of $\mathsf{Z}+\Pi_n\textsf{-Collection}$. \Square
\end{Theorems1}

\begin{Coroll1} \label{Th:LimitsOfCollectionInMinimalModels}
Let $n \geq 1$. Let $M$ be the minimal model of $\mathsf{Z}+\Pi_{n}\textsf{-Collection}$. Then there is an instance of $\Pi_{n+1}\textsf{-Collection}$ that fails in $\langle M, \in \rangle$.\Square
\end{Coroll1}

\begin{Theorems1}
Let $n \geq 1$. Let $M$ be the minimal model of $\mathsf{Z}+\Pi_{n}\textsf{-Collection}$. Then $\langle M, \in \rangle$ has a proper $\Sigma_{n+1}$-elementary end extension, but neither
\begin{itemize}
\item[(I)] a proper $\Sigma_{n+1}$-elementary end extension satisfying $\mathsf{KP}+\Pi_{n+3}\textsf{-Foundation}$, nor
\item[(II)] a proper $\Sigma_{n+1}$-elementary end extension satisfying $\mathsf{KP}+\Pi_n\textsf{-Collection}$.
\end{itemize}
\end{Theorems1}

\begin{proof}
The fact that $\langle M, \in \rangle$ has a proper $\Sigma_{n+1}$-elementary end extension follows from Theorem \ref{Th:KaufmannEEResult}. Let $\mathcal{N}= \langle N, \in^{\mathcal{N}} \rangle$ be such that $\mathcal{N} \models \mathsf{KP}$, $N \neq M$ and $\langle M, \in \rangle \prec_{e, n+1} \mathcal{N}$. Since $M$ is the minimal model of $\mathsf{Z}+\Pi_{n}\textsf{-Collection}$, $\langle M, \in \rangle \models \neg \sigma$ where $\sigma$ is the sentence
$$\exists x (x \textrm{ is transitive} \land \langle x, \in \rangle \models \mathsf{Z}+\Pi_{n}\textsf{-Collection}).$$
Since $\sigma$ is $\Sigma_1^{\mathsf{KP}}$ and $\langle M, \in \rangle \prec_{e, 1} \mathcal{N}$, $\mathcal{N} \models \neg \sigma$. Since $\mathcal{N} \models \mathsf{KP}$ and $M \neq N$, $\mathsf{Ord}^\mathcal{N} \backslash \mathsf{Ord}^{\langle M, \in\rangle}$ is nonempty. If $\gamma$ is the least element of $\mathsf{Ord}^\mathcal{N} \backslash \mathsf{Ord}^{\langle M, \in\rangle}$, then 
$$\mathcal{N} \models (\langle L_{\gamma}, \in \rangle \models \mathsf{Z}+\Pi_{n}\textsf{-Collection}),$$
which contradicts the fact that $\mathcal{N} \models \neg \sigma$. Therefore, $\mathsf{Ord}^\mathcal{N} \backslash \mathsf{Ord}^{\langle M, \in\rangle}$ is nonempty and contains no least element. Therefore, by Theorem \ref{Th:CollectionFromPartiallyElementaryEE1} and Corollary \ref{Th:LimitsOfCollectionInMinimalModels}, there must be both an instance of $\Pi_n\textsf{-Collection}$ and an instance of $\Pi_{n+3}\textsf{-Foundation}$ that fails in $\mathcal{N}$.     
\Square
\end{proof}

\section[Building partially-elementary end extensions]{Building partially-elementary end extensions} \label{sec:positiveresults}

In this section we use admissible covers and the Barwise Compactness Theorem to build partially-elementary end extensions that also satisfy any given recursive theory that holds in the base model. Barwise \cite{bar74} and \cite[Appendix]{bar75} introduces the machinery of admissible covers to apply infinitary compactness arguments to nonstandard countable models. The proof of \cite[Theorem A.4.1]{bar75} shows that for any countable model $\mathcal{M}$ of $\mathsf{KP}+\textsf{Foundation}$ and for any recursively enumerable $\mathcal{L}$-theory $T$ that holds in $\mathcal{M}$, $\mathcal{M}$ has proper end extension that satisfies $T$. By calibrating \cite[Appendix]{bar75}, Ressayre \cite[Theorem 2.15]{res87} shows that this result also holds for countable models of $\mathsf{KP}+\Sigma_1\textsf{-Foundation}$. 

\begin{Theorems1} \label{Th:RessayreEE1}
(Ressayre) Let $\mathcal{M}= \langle M, \in^\mathcal{M} \rangle$ be a countable model of $\mathsf{KP}+\Sigma_1\textsf{-Foundation}$. Let $T$ be a recursively enumerable theory such that $\mathcal{M} \models T$. Then there exists $\mathcal{N} \models T$ such that $\mathcal{M} \subseteq_e \mathcal{N}$ and $M \neq N$. \Square 
\end{Theorems1}

In \cite[2.17 Remarks]{res87}, Ressayre notes, without providing the details, that if $\mathcal{M}$ satisfies $\mathsf{KP}+\Pi_n\textsf{-Collection}+\Pi_{n+1}\cup \Sigma_{n+1}\textsf{-Foundation}$, then the end extension obtained in Theorem \ref{Th:RessayreEE1} can be guaranteed to be $\Sigma_n$-elementary. In this section, we work through the details of this result showing that the assumption that the model $\mathcal{M}$ being extended satisfies $\Pi_n\textsf{-Foundation}$ is not necessary. Our main result (Theorem \ref{Th:MainEndExtensionTheorem}) can be viewed as a generalisation of \cite[Theorem 5.3]{em22}, where admissible covers are used to build powerset-preserving end extension of countable models of set theory. Here we follow the presentation of admissible covers presented in \cite{em22}.  

In order to present admissible covers of (not necessarily well-founded) models of  extensions of $\mathsf{KP}$ we need to describe extensions of Kripke-Platek Set Theory that allow structures to appear as {\it urelements} in the domain of discourse. Let $\mathcal{L}^*$ be obtained from $\mathcal{L}$ by adding a new unary predicate $\mathsf{U}$, binary relation $\mathsf{E}$ and unary function symbol $\mathsf{F}$. Let $\mathcal{L}^*_\mathsf{S}$ be obtained from $\mathcal{L}^*$ by adding a new binary predicate $\mathsf{S}$. The intention is that $\mathsf{U}$ distinguishes objects that are urelements from objects that are sets, the urelements together with $\mathsf{E}$ form an $\mathcal{L}$-structure, and $\in$ is a membership relation between sets or urelments and sets. That is, the $\mathcal{L}^*$- and $\mathcal{L}^*_\mathsf{S}$-structures we will consider will be structures in the form $\mathfrak{A}_{\mathcal{M}}= \langle \mathcal{M}; A, \in^{\mathfrak{A}}, \mathsf{F}^\mathfrak{A}\rangle$ or $\mathfrak{A}_{\mathcal{M}}= \langle \mathcal{M}; A, \in^{\mathfrak{A}}, \mathsf{F}^\mathfrak{A}, \mathsf{S}^\mathfrak{A}\rangle$, where $\mathcal{M}= \langle M, \mathsf{E}^{\mathfrak{A}} \rangle$, $M$ is the extension of $\mathsf{U}$, $\mathsf{E}^{\mathfrak{A}} \subseteq M \times M$, $A$ is the extension of $\neg \mathsf{U}$ and $\in^\mathfrak{A} \subseteq (M \cup A) \times A$. 

Following \cite{bar75} we simplify the presentation of $\mathcal{L}^*$- and $\mathcal{L}^*_\mathsf{S}$-formulae by treating these languages as two-sorted instead of one-sorted and using the following conventions:
\begin{itemize}
\item The variables $p$, $q$, $r$, $p_1$, \ldots range over elements of the domain that satisfy $\mathsf{U}$;
\item the variables $a$, $b$, $c$, $a_1$, \ldots range over elements of the domain that satisfy $\neg \mathsf{U}$;
\item the variables $x$, $y$, $z$, $w$, $x_1$, \ldots range over all elements of the domain. 
\end{itemize}
So, $\forall a (\cdots )$ is an abbreviation of $\forall x(\neg \mathsf{U}(x) \Rightarrow \cdots)$, $\exists p(\cdots)$ is an abbreviation of $\exists x( \mathsf{U}(x) \land \cdots)$, etc. These conventions are used in the following $\mathcal{L}^*_\mathsf{S}$-axioms and -axiom schemes:
\begin{itemize}
\item[]({\sf Extensionality for sets}) $\forall a \forall b(a=b \iff \forall x(x \in a \iff x \in b))$
\item[]({\sf Pair}) $\forall x \forall y \exists a \forall z(z \in a \iff z= x \lor z=y)$
\item[]({\sf Union}) $\forall a \exists b (\forall y \in b)(\forall x \in y)(x \in b)$
\end{itemize}
Let $\Gamma$ be a class of $\mathcal{L}^*_\mathsf{S}$-formulae. 
\begin{itemize}
\item[]({\sf $\Gamma$-Separation}) For all $\phi(x, \vec{z})$ in $\Gamma$,
$$\forall \vec{z} \forall a \exists b \forall x(x \in b \iff (x \in a) \land \phi(x, \vec{z}))$$
\item[]({\sf $\Gamma$-Collection}) For all $\phi(x, y, \vec{z})$ in $\Gamma$,
$$\forall \vec{z} \forall a ((\forall x \in a) \exists y \phi(x, y, \vec{z}) \Rightarrow \exists b (\forall x \in a)(\exists y \in b) \phi(x, y, \vec{z}))$$
\item[]({\sf $\Gamma$-Foundation}) For all $\phi(x, \vec{z})$ in $\Gamma$,
$$\forall \vec{z}(\exists x \phi(x, \vec{z}) \Rightarrow \exists y(\phi(y, \vec{z}) \land (\forall w \in y) \neg \phi(w, \vec{z})))$$ 
\end{itemize}  
The interpretation of the function symbol $\mathsf{F}$ will map urelements, $p$, to sets, $a$, such that the $\mathsf{E}$-extension of $p$ is equal to the $\in$-extension of $a$. This is captured by the following axiom:
\begin{itemize}
\item[]($\dagger$) $\forall p \exists a(a= \mathsf{F}(p) \land \forall x( x \mathsf{E} p \iff x \in a)) \land \forall b(\mathsf{F}(b)= \emptyset)$
\end{itemize}
The following theory is the analogue of $\mathsf{KP}$ in the language $\mathcal{L}^*$:
\begin{itemize}
\item $\mathsf{KPU}_{\mathbb{C}ov}$ is the $\mathcal{L}^*$-theory with axioms: $\exists a (a=a)$, $\forall p \forall x(x \notin p)$, {\sf Extensionality for sets}, {\sf Pair}, {\sf Union}, $\Delta_0(\mathcal{L}^*)\textsf{-Separation}$, $\Delta_0(\mathcal{L}^*)\textsf{-Collection}$, $\Pi_1(\mathcal{L}^*)\textsf{-Foundation}$ and ($\dagger$).
\end{itemize}
An order pair $\langle x, y \rangle$ is coded in $\mathsf{KPU}_{\mathbb{C}ov}$ by the set $\{\{x\}, \{x, y\}\}$, and we write $\mathsf{OP}(x)$ for the usual $\Delta_0$-formula that says that $z$ is an order pair and that also works in this theory. We write $\mathsf{fst}$ for the function $\langle x, y \rangle \mapsto x$ and $\mathsf{snd}$ for the function $\langle x, y \rangle \mapsto y$. The usual $\Delta_0$ definitions of the graphs of these functions also work in $\mathsf{KPU}_{\mathbb{C}ov}$. The {\bf rank function}, $\rho$, and {\bf support function}, $\mathsf{sp}$, are defined in $\mathsf{KPU}_{\mathbb{C}ov}$ by recursion:
$$\rho(p)= 0 \textrm{ for all urelements } p, \textrm{ and } \rho(a)= \sup \{\rho(x)+1 \mid x \in a\} \textrm{ for all sets } a;$$
$$\mathsf{sp}(p)=\{p\} \textrm{ for all urelements } p, \textrm{ and } \mathsf{sp}(a)= \bigcup_{x \in a} \mathsf{sp}(x) \textrm{ for all sets } a.$$
The theory $\mathsf{KPU}_{\mathbb{C}ov}$ proves that both $\mathsf{sp}$ and $\rho$ are total functions and their graphs are $\Delta_1(\mathcal{L}^*)$. We say that $x$ is a {\bf pure set} if $\mathsf{sp}(x)=\emptyset$. The following $\Delta_0(\mathcal{L}^*)$-formulae assert that `$x$ is transitive' and `$x$ is an ordinal (a hereditarily transitive pure set)':
$$\mathsf{Transitive}(x) \iff \neg\mathsf{U}(x) \land (\forall y \in x)(\forall z \in y) (z \in x);$$
$$\mathsf{Ord}(x) \iff (\mathsf{Transitive}(x) \land (\forall y \in x) \mathsf{Transitive}(y))$$
We will consider $\mathcal{L}^*_\mathsf{S}$-structures in which the predicate $\mathsf{S}$ is a satisfaction class for the $\Sigma_n$-formulae of the $\mathcal{L}$-structure $\mathcal{M}$. Let $\mathsf{KPU}^\prime_{\mathbb{C}ov}$ be obtained from $\mathsf{KPU}_{\mathbb{C}ov}$ by adding axioms asserting that the $\mathcal{L}$-structure formed by the urelements and the binary relation $\mathsf{E}$ satisfies $\mathsf{KP}$. For $n \in \omega$, define
\begin{itemize}
\item[]($n\textsf{-Sat}$) $\mathsf{S}(m, x)$ if and only if $\mathsf{U}(m)$ and $\mathsf{U}(x)$ and $\mathsf{Sat}_{\Sigma_n}(m, x)$ holds in the $\mathcal{L}$-structure defined by $\mathsf{U}$ and $\mathsf{E}$.
\end{itemize}
We can now define a family of $\mathcal{L}^*_\mathsf{S}$-theories extending $\mathsf{KPU}_{\mathbb{C}ov}$ that assert that the $\mathcal{L}$-structure defined by $\mathsf{U}$ and $\mathsf{E}$ satisfies $\mathsf{KP}$ and $\mathsf{S}$ is a satisfaction class on this structure for $\Sigma_n$-formulae, and $\mathsf{S}$ can be used in the separation, collection and foundation schemes. 
\begin{itemize}
\item For all $n \in \omega$, define $\mathsf{KPU}^n_{\mathbb{C}ov}$ to be the $\mathcal{L}^*_\mathsf{S}$-theory extending $\mathsf{KPU}^\prime_{\mathbb{C}ov}$ with the axiom $n\textsf{-Sat}$ and the schemes $\Delta_0(\mathcal{L}^*_\mathsf{S})\textsf{-Separation}$, $\Delta_0(\mathcal{L}^*_\mathsf{S})\textsf{-Collection}$ and $\Pi_1(\mathcal{L}^*_\mathsf{S})\textsf{-Foundation}$.
\end{itemize}
The arguments used in \cite[I.4.4 and I.4.5]{bar75} show that $\mathsf{KPU}_{\mathbb{C}ov}$ proves the schemes of $\Sigma_1(\mathcal{L}^*)\textsf{-Collection}$ and $\Delta_1(\mathcal{L}^*)\textsf{-Separation}$, and for all $n \in \omega$, $\mathsf{KPU}^n_{\mathbb{C}ov}$ proves the schemes of $\Sigma_1(\mathcal{L}^*_\mathsf{S})\textsf{-Collection}$ and $\Delta_1(\mathcal{L}^*_\mathsf{S})\textsf{-Separation}$. 

\begin{Definitions1}
Let $\mathcal{M}= \langle M, \mathsf{E}^\mathcal{M} \rangle$ be an $\mathcal{L}$-structure. An {\bf admissible set covering $\mathcal{M}$} is an $\mathcal{L}^*$-structure
$$\mathfrak{A}_\mathcal{M} = \langle \mathcal{M}; A, \in^\mathfrak{A}, \mathsf{F}^\mathfrak{A} \rangle \models \mathsf{KPU}_{\mathbb{C}ov}$$
such that $\in^\mathfrak{A}$ is well-founded. If $\mathcal{M} \models \mathsf{KP}$ and $n \in \omega$, then an {\bf $n$-admissible set covering $\mathcal{M}$} is an $\mathcal{L}^*_\mathsf{S}$-structure
$$\mathfrak{A}_\mathcal{M} = \langle \mathcal{M}; A, \in^\mathfrak{A}, \mathsf{F}^\mathfrak{A}, \mathsf{S}^\mathfrak{A} \rangle \models \mathsf{KPU}^n_{\mathbb{C}ov}$$
such that $\in^\mathfrak{A}$ is well-founded. Note that if $\mathfrak{A}_\mathcal{M} = \langle \mathcal{M}; A, \in^\mathfrak{A}, \mathsf{F}^\mathfrak{A}, \ldots \rangle$ is an ($n$-)admissible set covering $\mathcal{M}$, then $\mathfrak{A}_\mathcal{M}$ is isomorphic to a structure whose membership relation ($\in$) is the membership relation of the metatheory. The {\bf admissible cover} of $\mathcal{M}$, denoted $\mathbb{C}\mathrm{ov}_\mathcal{M}= \langle \mathcal{M}; A_\mathcal{M}, \in, \mathsf{F}_\mathcal{M} \rangle$, is the smallest admissible set covering $\mathcal{M}$ whose membership relation ($\in$) coincides with the membership relation of the metatheory. If $\mathcal{M} \models \mathsf{KP}$ and $n \in \omega$, the {\bf $n$-admissible cover} of $\mathcal{M}$, denoted $\mathbb{C}\mathrm{ov}^n_\mathcal{M}= \langle \mathcal{M}; A_\mathcal{M}, \in, \mathsf{F}_\mathcal{M}, \mathsf{S}_\mathcal{M} \rangle$, is the smallest $n$-admissible set covering $\mathcal{M}$ whose membership relation ($\in$) coincides with the membership relation of the metatheory.  
\end{Definitions1}

\begin{Definitions1}
Let $\mathcal{M}= \langle M, \mathsf{E}^\mathcal{M} \rangle$ be an $\mathcal{L}$-structure, and let $\mathfrak{A}_\mathcal{M} = \langle \mathcal{M}; A, \in, \mathsf{F}^\mathfrak{A}, \ldots \rangle$ be an $\mathcal{L}^*$- or $\mathcal{L}^*_\mathsf{S}$-structure. We use $\mathrm{WF}(A)$ to denote the largest $B \subseteq A$ such that $\langle B, \in^\mathfrak{A} \rangle \subseteq_e \langle A, \in^\mathfrak{A} \rangle$ and $\langle B, \in^\mathfrak{A} \rangle$ is well-founded. The {\bf well-founded part} of $\mathfrak{A}_\mathcal{M}$ is the $\mathcal{L}^*$- or $\mathcal{L}^*_\mathsf{S}$-structure
$$\mathrm{WF}(\mathfrak{A}_\mathcal{M})= \langle \mathcal{M}; \mathrm{WF}(A), \in^\mathfrak{A}, \mathsf{F}^\mathfrak{A}, \ldots \rangle$$
Note that $\mathrm{WF}(\mathfrak{A}_\mathcal{M})$ is always isomorphic to a structure whose membership relation $\in$ coincides with the membership relation of the metatheory.     
\end{Definitions1}

Let $\mathcal{M}= \langle M, \mathsf{E}^\mathcal{M} \rangle$ be such that $\mathcal{M} \models \mathsf{KP}$. Let $\mathcal{L}^\mathtt{ee}$ be the language obtained from $\mathcal{L}$ by adding new constant symbols $\bar{a}$ for each $a \in M$ and a new constant symbol $\mathbf{c}$. Let $\mathfrak{A}_\mathcal{M}= \langle \mathcal{M}; A, \in, \mathsf{F}^\mathfrak{A}, \mathsf{S}^\mathfrak{A} \rangle$ be an $n$-admissible set covering $\mathcal{M}$. There is a coding $\ulcorner \cdot \urcorner$ of a fragment of the infinitary language $\mathcal{L}_{\infty \omega}^\mathtt{ee}$ in $\mathfrak{A}_\mathcal{M}$ with the property that the classes of codes of {\it atomic formulae}, {\it variables}, {\it constants}, {\it well-formed formulae}, {\it sentences}, etc. are all $\Delta_1(\mathcal{L}^*)$-definable over $\mathfrak{A}_\mathcal{M}$ (see \cite[p9]{em22} for an explicit definition of such a coding). We write $\mathcal{L}_{\mathfrak{A}_\mathcal{M}}^\mathtt{ee}$ for the fragment of $\mathcal{L}_{\infty \omega}^\mathtt{ee}$ whose codes appear in $\mathfrak{A}_\mathcal{M}$. In order to apply compactness arguments to $\mathcal{L}_{\mathfrak{A}_\mathcal{M}}^\mathtt{ee}$-theories where $\mathfrak{A}_\mathcal{M}$ is an $n$-admissible set, we will use the following specific version of the Barwise Compactness Theorem (\cite[III.5.6]{bar75}):

\begin{Theorems1}
(Barwise Compactness Theorem) Let $\mathfrak{A}_\mathcal{M}= \langle \mathcal{M}; A, \in \mathsf{F}^\mathfrak{A}, \mathsf{S}^\mathfrak{A} \rangle$ be an $n$-admissible set covering $\mathcal{M}$. Let $T$ be an $\mathcal{L}_{\mathfrak{A}_\mathcal{M}}^\mathtt{ee}$-theory that is $\Sigma_1(\mathcal{L}^*_\mathsf{S})$-definable over $\mathfrak{A}_\mathcal{M}$ and such that for all $T_0 \subseteq T$, if $T_0 \in A$, then $T_0$ has a model. Then $T$ has a model.\Square  
\end{Theorems1}

The work in \cite[Appendix]{bar75} and \cite[Chapter 2]{res87} shows that if $\mathcal{M}$ satisfies $\mathsf{KP}+\Sigma_1\textsf{-Foundation}$, then $\mathbb{C}\mathrm{ov}_\mathcal{M}$ exists. In particular, $\mathbb{C}\mathrm{ov}_\mathcal{M}$ can be obtained from $\mathcal{M}$ by first defining a model of $\mathsf{KPU}_{\mathbb{C}ov}$ inside $\mathcal{M}$ and then considering the well-founded part of this model. We now turn to reviewing the construction of $\mathbb{C}\mathrm{ov}_\mathcal{M}$ from $\mathcal{M}$ and showing that if $\mathcal{M}$ satisfies $\mathsf{KP}+\Pi_n\textsf{-Collection}+\Sigma_{n+1}\textsf{-Foundation}$, then $\mathbb{C}\mathrm{ov}_\mathcal{M}$ can be expanded to an $\mathcal{L}^*_\mathsf{S}$-structure corresponding to $\mathbb{C}\mathrm{ov}^n_\mathcal{M}$.

Let $n \geq 1$. Fix a model $\mathcal{M}= \langle M, \mathsf{E}^\mathcal{M} \rangle$ that satisfies $\mathsf{KP}+\Pi_n\textsf{-Collection}+\Sigma_{n+1}\textsf{-Foundation}$. Working inside $\mathcal{M}$, define unary relations $\mathsf{N}$ and $\mathsf{Set}$, binary relations $\mathsf{E}^\prime$, $\mathcal{E}$ and $\bar{\mathsf{S}}$, and unary function $\bar{F}$ by:
$$\mathsf{N}(x) \textrm{ iff } \exists y(x = \langle 0, y \rangle);$$
$$x \mathsf{E}^\prime y \textrm{ iff } \exists w \exists z(x= \langle 0, w \rangle \land y= \langle 0, z \rangle \land w \in z);$$
$$\mathsf{Set}(x)= \exists y (x = \langle 1, y \rangle \land (\forall z \in y) (\mathsf{N}(z) \lor \mathsf{Set}(z)));$$
$$x \mathcal{E} y \textrm{ iff } \exists z (y= \langle 1, z \rangle \land x \in z);$$
$$\bar{\mathsf{F}}(x)= \langle 1, X \rangle \textrm{ where } X= \{ \langle 0, y \rangle \mid \exists w (x= \langle 0, w \rangle \land y \in w)\};$$
$$\bar{\mathsf{S}}(x, y) \textrm{ iff } \exists z \exists w(x= \langle 0, w \rangle \land y= \langle 0, z \rangle \land \mathsf{Sat}_{\Sigma_n}(w, z)).$$
It is noted in \cite[Appendix Section 3]{bar75} that $\mathsf{N}$, $\mathsf{E}^\prime$, $\mathcal{E}$ and $\bar{\mathsf{F}}$ are defined by $\Delta_0$-formulae in $\mathcal{M}$. The Second Recursion Theorem (\cite[V.2.3.]{bar75}), provable in $\mathsf{KP}+\Sigma_1\textsf{-Foundation}$ as note in \cite{res87}, ensures that $\mathsf{Set}$ can be expressed as a $\Sigma_1$-formula in $\mathcal{M}$. Theorem \ref{Complexityofpartialsat} implies that $\bar{\mathsf{S}}$ is defined by a $\Sigma_n$-formula in $\mathcal{M}$. These definitions yield an interpretation, $\mathcal{I}$, of an $\mathcal{L}^*_\mathsf{S}$-structure $\mathfrak{A}_\mathcal{N}= \langle \mathcal{N}; \mathsf{Set}^\mathcal{M}, \mathcal{E}^\mathcal{M}, \bar{\mathsf{F}}^\mathcal{M}, \bar{\mathsf{S}}^\mathcal{M} \rangle$, where $\mathcal{N}= \langle \mathsf{N}^\mathcal{M}, (\mathsf{E}^\prime)^\mathcal{M} \rangle$. The following table that extends the table on page 373 in \cite{bar75} summarises the interpretation $\mathcal{I}$:
\begin{table}[H]
  \begin{center}
    \caption{The interpretation $\mathcal{I}$}
    \label{tab:table1}
    \begin{tabular}{c|c} % <-- Alignments: 1st column left, 2nd middle and 3rd right, with vertical lines in between
      {\bf $\mathcal{L}^*_\mathsf{S}$ Symbol} & {\bf $\mathcal{L}$ expression under $\mathcal{I}$}\\
      \hline\\
      $\forall x$ & $\forall x(\mathsf{N}(x) \lor \mathsf{Set}(x) \Rightarrow \cdots)$ \\
      $=$ & $=$\\
      $\mathsf{U}(x)$ & $\mathsf{N}(x)$\\
      %$S(x)$ & $\mathsf{Set}(x)$\\
      $x \mathsf{E} y$ & $x \mathsf{E}^\prime y$\\
      $x \in y$ & $x \mathcal{E} y$\\
      $\mathsf{F}(x)$ & $\bar{\mathsf{F}}(x)$\\
	  $\mathsf{S}(x, y)$ & $\bar{\mathsf{S}}(x, y)$
    \end{tabular}
  \end{center}
\end{table}
If $\phi$ is an $\mathcal{L}^*_\mathsf{S}$-formula, then we write $\phi^\mathcal{I}$ for the translation of $\phi$ into an $\mathcal{L}$-formula described in this table. By ignoring the interpretation $\bar{\mathsf{S}}$ of $\mathsf{S}$ we obtain, instead, an interpretation, $\mathcal{I}^-$, of an $\mathcal{L}^*$-structure in $\mathcal{M}$ and we write $\mathfrak{A}_\mathcal{N}^-$ for this reduct. Note that the map $x \mapsto \langle 0, x \rangle$ defines an isomorphism between $\mathcal{M}$ and $\mathcal{N}= \langle \mathsf{N}^\mathcal{M}, (\mathsf{E}^\prime)^\mathcal{M} \rangle$. Ressayre, refining \cite[Appendix Lemma 3.2]{bar75}, shows that if $\mathcal{M}$ satisfies $\mathsf{KP}+\Sigma_1\textsf{-Foundation}$, then interpretation $\mathcal{I}^-$ yields a structure satisfying $\mathsf{KPU}_{\mathbb{C}\mathrm{ov}}$.

\begin{Theorems1}
$\mathfrak{A}_\mathcal{N}^- \models \mathsf{KPU}_{\mathbb{C}\mathrm{ov}}$. \Square
\end{Theorems1}

\begin{Lemma1}
Let $\phi(\vec{x})$ be a $\Delta_0(\mathcal{L}_\mathsf{S}^*)$-formula. Then $\phi^{\mathcal{I}}(\vec{x})$ is equivalent to a $\Delta_{n+1}$-formula in $\mathcal{M}$.
\end{Lemma1}

\begin{proof}
We prove this result by induction on the complexity of $\phi$. Above, we observed that $\mathsf{N}(x)$, $x\mathsf{E}^\prime y$, $x \mathcal{E} y$ and $y= \bar{\mathsf{F}}(x)$ can be written as $\Delta_0$-formulae. And $\bar{\mathsf{S}}(x, y)$ can be written as a $\Sigma_n$-formula. Now, $y \mathcal{E} \bar{F}(x)$ if and only if
$$\mathsf{fst}(y)=0 \land \mathsf{snd}(y) \in \mathsf{snd}(x),$$
which is $\Delta_0$. Therefore, if $\phi(\vec{x})$ is a quantifier-free $\mathcal{L}_\mathsf{S}^*$-formula, then $\phi^{\mathcal{I}}(\vec{x})$ is equivalent to a $\Delta_{n+1}$-formula in $\mathcal{M}$. Now, suppose that $\phi(x_0, \ldots, x_{m-1})$ is in the form $(\exists y \in x_0) \psi(x_0, \ldots, x_{m-1}, y)$ where $\psi^{\mathcal{I}}(x_0, \ldots, x_{m-1}, y)$ is equivalent to a $\Delta_{n+1}$-formula in $\mathcal{M}$. Therefore, $\phi^{\mathcal{I}}(x_0, \ldots, x_{m-1})= (\exists y \mathcal{E} x_0) \psi^{\mathcal{I}}(x_0, \ldots, x_{m-1}, y)$, and $(\exists y \mathcal{E} x_0) \psi^{\mathcal{I}}(x_0, \ldots, x_{m-1}, y)$ iff
$$(\exists y \in \mathsf{snd}(x_0)) \psi^{\mathcal{I}}(x_0, \ldots, x_{m-1}, y)$$
So, since $\mathcal{M}$ satisfies $\Pi_n\textsf{-Collection}$, $\phi^{\mathcal{I}}(x_0, \ldots, x_{m-1})$ is equivalent to a $\Delta_{n+1}$-formula in $\mathcal{M}$. Finally, suppose that $\phi(x_0, \ldots, x_{m-1})$ is in the form $(\exists y \in \mathsf{F}(x_0)) \psi(x_0, \ldots, x_{m-1}, y)$ where $\psi^{\mathcal{I}}(x_0, \ldots, x_{m-1}, y)$ is equivalent to a $\Delta_{n+1}$-formula in $\mathcal{M}$. Therefore, $\phi^{\mathcal{I}}(x_0, \ldots, x_{m-1})= (\exists y \mathcal{E} \bar{\mathsf{F}}(x_0)) \psi^{\mathcal{I}}(x_0, \ldots, x_{m-1}, y)$, and $(\exists y \mathcal{E} \bar{\mathsf{F}}(x_0)) \psi^{\mathcal{I}}(x_0, \ldots, x_{m-1}, y)$ iff
$$\exists z (z= \bar{\mathsf{F}}(x_0) \land (\exists y \in \mathsf{snd}(z)) \psi^{\mathcal{I}}(x_0, \ldots, x_{m-1}, y))$$
$$\textrm{iff } \forall z(z= \bar{\mathsf{F}}(x_0) \Rightarrow (\exists y \in \mathsf{snd}(z)) \psi^{\mathcal{I}}(x_0, \ldots, x_{m-1}, y))$$
Therefore, since $\mathcal{M}$ satisfies $\Pi_n$-Collection, $\phi^{\mathcal{I}}(x_0, \ldots, x_{m-1})$ is equivalent to a $\Delta_{n+1}$-formula in $\mathcal{M}$. The Lemma now follows by induction. \Square       
\end{proof}

\begin{Lemma1}
$\mathfrak{A}_\mathcal{N} \models \Delta_0(\mathcal{L}^*_{\mathsf{S}})\textsf{-Separation}$.
\end{Lemma1}

\begin{proof}
Let $\phi(x, \vec{z})$ be a $ \Delta_0(\mathcal{L}^*_{\mathsf{S}})$-formula. Let $\vec{v}$ be sets and/or urelements of $\mathfrak{A}_\mathcal{N}$ and $a$ a set of $\mathfrak{A}_\mathcal{N}$. Work inside $\mathcal{M}$. Now, $a= \langle 1, a_0 \rangle$. Let
\[
b_0=\{x \in a_0 \mid \phi^{\mathcal{I}}(x, \vec{v})\},
\] 
which is a set by $\Delta_{n+1}\textsf{-Separation}$. Let $b=\langle 1, b_0 \rangle$. Therefore, for all $x$ such that $\mathsf{Set}(x)$,
\[
x \mathcal{E} b \textrm{ if and only if } x \mathcal{E} a \land \phi^{\mathcal{I}}(x, \vec{v}).
\]
This shows that $\mathfrak{A}_\mathcal{N}$ satisfies $\Delta_0(\mathcal{L}^*_{\mathsf{S}})\textsf{-Separation}$.\Square 
\end{proof}

\begin{Lemma1}
$\mathfrak{A}_\mathcal{N} \models \Delta_0(\mathcal{L}^*_{\mathsf{S}})\textsf{-Collection}$.
\end{Lemma1}
       
\begin{proof}
Let $\phi(x, y, \vec{z})$ be a $\Delta_0(\mathcal{L}^*_{\mathsf{S}})$-formula. Let $\vec{v}$ be a sequence of sets and/or urelements of $\mathfrak{A}_{\mathcal{N}}$ and let $a$ be a set of $\mathfrak{A}_{\mathcal{N}}$ such that 
\[
\mathfrak{A}_{\mathcal{N}} \models (\forall x \in a) \exists y \phi(x, y, \vec{v}).
\]
Work inside $\mathcal{M}$. Now, $a= \langle 1, a_0 \rangle$. And,
\[
(\forall x \mathcal{E} a) \exists y ((\mathsf{N}(y) \lor \mathsf{Set}(y)) \land \phi^{\mathcal{I}}(x, y, \vec{v})).
\] 
So,
\[
(\forall x \in a_0) \exists y((\mathsf{N}(y) \lor \mathsf{Set}(y)) \land \phi^{\mathcal{I}}(x, y, \vec{v})).
\]
Since $(\mathsf{N}(y) \lor \mathsf{Set}(y)) \land \phi^{\mathcal{I}}(x, y, \vec{v})$ is equivalent to a $\Sigma_{n+1}$-formula, we can use $\Pi_n\textsf{-Collection}$ to find $b_0$ such that
\[
(\forall x \in a_0) (\exists y \in b_0)((\mathsf{N}(y) \lor \mathsf{Set}(y)) \land \phi^{\mathcal{I}}(x, y, \vec{v})).
\]
Let $b_1= \{y \in b_0 \mid \mathsf{N}(y) \lor \mathsf{Set}(y)\}$, which is a set by $\Sigma_1\textsf{-Separation}$. Let $b= \langle 1, b_1 \rangle$. Therefore, $\mathsf{Set}(b)$ and 
\[
(\forall x \mathcal{E} a)(\exists y \mathcal{E} b) \phi^{\mathcal{I}}(x, y, \vec{v}).
\] 
So, 
\[
\mathfrak{A}_{\mathcal{N}} \models (\forall x \in a) (\exists y \in b) \phi(x, y, \vec{v}).
\]
This shows that $\mathfrak{A}_\mathcal{N}$ satisfies $\Delta_0(\mathcal{L}^*_{\mathsf{S}})\textsf{-Collection}$.
\Square
\end{proof}       

\begin{Lemma1}
$\mathfrak{A}_\mathcal{N} \models \Sigma_1(\mathcal{L}^*_{\mathsf{S}})\textsf{-Foundation}$.
\end{Lemma1}

\begin{proof}
Let $\phi(x, \vec{z})$ be a $\Sigma_1(\mathcal{L}^*_{\mathsf{S}})$-formula. Let $\vec{v}$ be a sequence of sets and/or urelements such that
\[
\{x \in \mathfrak{A}_{\mathcal{N}} \mid \mathfrak{A}_\mathcal{N} \models \phi(x, \vec{v}) \} \textrm{ is nonempty}.
\]
Work inside $\mathcal{M}$. Consider $\theta(\alpha, \vec{z})$ defined by
\[
(\alpha \textrm{ is an ordinal}) \land \exists x((\mathsf{N}(x) \lor \mathsf{Set}(x)) \land \rho(x)= \alpha \land \phi^{\mathcal{I}}(x, \vec{z})).
\]
Note that $\theta(\alpha, \vec{z})$ is equivalent to a $\Sigma_{n+1}$-formula and $\exists \alpha \theta(\alpha, \vec{v})$. Therefore, using $\Sigma_{n+1}\textsf{-Foundation}$, let $\beta$ be a $\in$-least element of 
\[
\{ \alpha \in M \mid \mathcal{M} \models \theta(\alpha, \vec{v})\}.
\]
Let $y$ be such that $(\mathsf{N}(y) \lor \mathsf{Set}(y))$, $\rho(y)=\beta$ and $\phi^{\mathcal{I}}(y, \vec{v})$. Note that if $x \mathcal{E} y$, then $\rho(x) < \rho(y)$. Therefore $y$ is an $\mathcal{E}$-least element of
\[
\{x \in \mathfrak{A}_{\mathcal{N}} \mid \mathfrak{A}_\mathcal{N} \models \phi(x, \vec{v}) \}.
\] 
\Square 
\end{proof}

The results of \cite[Appendix Section 3]{bar75} show that $\mathbb{C}\mathrm{ov}_\mathcal{M}$ is the $\mathcal{L}^*$-reduct of the well-founded part of $\mathfrak{A}_\mathcal{N}$.

\begin{Theorems1} \label{Th:CharacterisationOfCov}
(Barwise) The $\mathcal{L}^*$-reduct of $\mathrm{WF}(\mathfrak{A}_\mathcal{N})$, $\mathrm{WF}^-(\mathfrak{A}_\mathcal{N})= \langle \mathcal{N}; \mathrm{WF}(\mathsf{Set}^\mathcal{M}), \mathcal{E}^\mathcal{M}, \bar{\mathsf{F}}^\mathcal{M} \rangle$, is an admissible set covering $\mathcal{N}$ that is isomorphic to $\mathbb{C}\mathrm{ov}_\mathcal{M}$. \Square
\end{Theorems1}

We can extend this result to show that $\mathrm{WF}(\mathfrak{A}_\mathcal{N})$ is an $n$-admissible cover of $\mathcal{N}$ and, therefore, isomorphic to $\mathbb{C}\mathrm{ov}^n_\mathcal{M}$.

\begin{Theorems1} 
The structure $\mathrm{WF}(\mathfrak{A}_\mathcal{N})= \langle \mathcal{N}; \mathrm{WF}(\mathsf{Set}^\mathcal{M}), \mathcal{E}^\mathcal{M}, \bar{\mathsf{F}}^\mathcal{M}, \bar{\mathsf{S}}^\mathcal{M} \rangle$ is an $n$-admissible set covering $\mathcal{N}$. Moreover, $\mathrm{WF}(\mathfrak{A}_\mathcal{N})$ is isomorphic to $\mathbb{C}\mathrm{ov}^n_\mathcal{M}$.
\end{Theorems1}

\begin{proof}
Theorem \ref{Th:CharacterisationOfCov}, the fact that $\mathcal{M} \models \mathsf{KP}$, and the fact that $\mathrm{WF}(\mathfrak{A}_\mathcal{N})$ is well-founded imply that $\mathrm{WF}(\mathfrak{A}_\mathcal{N})$ satisfies $\mathsf{KPU}^\prime_{\mathbb{C}\mathrm{ov}}+\mathcal{L}^*_\mathsf{S}\textsf{-Foundation}$. The definition of $\bar{\mathsf{S}}$ in $\mathcal{M}$ ensures that $\mathrm{WF}(\mathfrak{A}_\mathcal{N})$ satisfies $n\textsf{-Sat}$. If $a$ is a set $\mathrm{WF}(\mathfrak{A}_\mathcal{N})$ and $b$ is a set in $\mathfrak{A}_\mathcal{N}$ with $\mathfrak{A}_\mathcal{N} \models (b \subseteq a)$, then $b \in \mathrm{WF}(\mathsf{Set}^\mathcal{M})$. Therefore, since $\Delta_0(\mathcal{L}^*_{\mathsf{S}})$-formulae are absolute between $\mathrm{WF}(\mathfrak{A}_\mathcal{N})$ and $\mathfrak{A}_\mathcal{N}$, $\mathrm{WF}(\mathfrak{A}_\mathcal{N})$ satisfies $\Delta_0(\mathcal{L}^*_{\mathsf{S}})\textsf{-Separation}$. To show that $\mathrm{WF}(\mathfrak{A}_\mathcal{N})$ satisfies $\Delta_0(\mathcal{L}^*_{\mathsf{S}})\textsf{-Collection}$, let $\phi(x, y, \vec{z})$ be a $\Delta_0(\mathcal{L}^*_{\mathsf{S}})$-formula. Let $\vec{v}$ be sets and/or urelements in $\mathrm{WF}(\mathfrak{A}_\mathcal{N})$ and let $a$ be a set of $\mathrm{WF}(\mathfrak{A}_\mathcal{N})$ such that
\[
\mathrm{WF}(\mathfrak{A}_\mathcal{N}) \models (\forall x \in a)\exists y \phi(x, y, \vec{v}).
\] 
Consider the formula $\theta(\beta, \vec{z})$ defined by
\[
(\beta \textrm{ is an ordinal}) \land (\forall x \in a)(\exists \alpha \in \beta)\exists y (\rho(y)=\alpha \land \phi(x, y, \vec{z})).
\]
Note that if $\beta$ is a nonstandard ordinal of $\mathfrak{A}_\mathcal{N}$, then $\mathfrak{A}_\mathcal{N} \models \theta(\beta, \vec{v})$. Using $\Delta_0(\mathcal{L}^*_{\mathsf{S}})\textsf{-Collection}$, $\theta(\beta, \vec{z})$ is equivalent to a $\Sigma_1(\mathcal{L}^*_{\mathsf{S}})$-formula in $\mathfrak{A}_\mathcal{N}$. Therefore, by $\Sigma_1(\mathcal{L}^*_{\mathsf{S}})\textsf{-Foundation}$ in $\mathfrak{A}_\mathcal{N}$, $\{ \beta \mid \mathfrak{A}_\mathcal{N} \models \theta(\beta, \vec{v})\}$ has a least element $\gamma$. Note that $\gamma$ must be an ordinal in $\mathrm{WF}(\mathfrak{A}_\mathcal{N})$. Consider the formula $\psi(x, y, \vec{z}, \gamma)$ defined by $\phi(x, y, \vec{z}) \land (\rho(y) < \gamma)$. Then,
\[
\mathfrak{A}_\mathcal{N} \models (\forall x \in a) \exists y \psi(x, y, \vec{v}, \gamma).
\]
By $\Delta_0(\mathcal{L}^*_{\mathsf{S}})\textsf{-Collection}$ in $\mathfrak{A}_\mathcal{N}$, there is a set $b$ of $\mathfrak{A}_\mathcal{N}$ such that
\[
\mathfrak{A}_\mathcal{N} \models (\forall x \in a) (\exists y \in b) \psi(x, y, \vec{v}, \gamma).
\]
Let $c= \{ y \in b \mid \rho(y)< \gamma\}$, which is a set in $\mathfrak{A}_\mathcal{N}$ by $\Delta_1(\mathcal{L}^*_{\mathsf{S}})\textsf{-Separation}$. Now, $c$ is a set of $\mathrm{WF}(\mathfrak{A}_\mathcal{N})$ and
\[
\mathrm{WF}(\mathfrak{A}_\mathcal{N}) \models (\forall x \in a) (\exists y \in c) \phi(x, y, \vec{v}).
\]
Therefore, $\mathrm{WF}(\mathfrak{A}_\mathcal{N})$ satisfies $\Delta_0(\mathcal{L}^*_{\mathsf{S}})\textsf{-Collection}$, and so is an $n$-admissible set covering $\mathcal{N}$. Since the $\mathcal{L}^*$-reduct of $\mathrm{WF}(\mathfrak{A}_\mathcal{N})$ is isomorphic to $\mathbb{C}\mathrm{ov}_\mathcal{M}$, $\mathrm{WF}(\mathfrak{A}_\mathcal{N})$ is isomorphic to $\mathbb{C}\mathrm{ov}^n_\mathcal{M}$.
\Square       
\end{proof}

To summarise, we have proved the following:

\begin{Theorems1}
If $\mathcal{M}\models \mathsf{KP}+\Pi_n\textsf{-Collection}+\Sigma_{n+1}\textsf{-Foundation}$, then then there is an interpretation of $\mathsf{S}$ in $\mathbb{C}\mathrm{ov}_\mathcal{M}$ that yields the $n$-admissible cover $\mathbb{C}\mathrm{ov}^n_\mathcal{M}$.\Square   
\end{Theorems1}

Our analysis also yields the following version of \cite[Appendix Corollary 2.4]{bar75}, which plays an important role on compactness arguments:

\begin{Theorems1}
Let $\mathcal{M}= \langle M, \mathsf{E}^\mathcal{M} \rangle$ be such that $\mathcal{M}\models \mathsf{KP}+\Pi_n\textsf{-Collection}+\Sigma_{n+1}\textsf{-Foundation}$. For all $A \subseteq M$, there exists $a \in M$ such that $a^*=A$ if and only if $A \in \mathbb{C}\mathrm{ov}^n_\mathcal{M}$.\Square
\end{Theorems1}

In particular, we obtain:

\begin{Lemma1} \label{Th:KeyCompactnessLemma}
Let $\mathcal{M}= \langle M, \mathsf{E}^\mathcal{M} \rangle$ be such that $\mathcal{M}\models \mathsf{KP}+\Pi_n\textsf{-Collection}+\Sigma_{n+1}\textsf{-Foundation}$. Let $T_0$ be an $\mathcal{L}_{\mathbb{C}\mathrm{ov}^n_\mathcal{M}}^{\mathtt{ee}}$-theory. If $T_0 \in \mathbb{C}\mathrm{ov}^n_\mathcal{M}$, then there exists $b \in M$ such that
\[
b^*= \{ a \in M \mid \bar{a} \textrm{ is mentioned in } T_0\}.
\]
\Square
\end{Lemma1}

The next result connects definability in $\mathcal{M}$ with definability in $\mathbb{C}\mathrm{ov}^n_\mathcal{M}$.

\begin{Lemma1}
Let $\mathcal{M}= \langle M, \mathsf{E}^\mathcal{M} \rangle$ be such that $\mathcal{M}\models \mathsf{KP}+\Pi_n\textsf{-Collection}+\Sigma_{n+1}\textsf{-Foundation}$. Let $\phi(\vec{z})$ be a $\Sigma_{n+1}$-formula. Then there exists a $\Sigma_1(\mathcal{L}^*_{\mathsf{S}})$-formula $\hat{\phi}(\vec{z})$ such that for all $\vec{z} \in M$,
\[
\mathcal{M} \models \phi(\vec{z}) \textrm{ if and only if } \mathbb{C}\mathrm{ov}^n_\mathcal{M} \models \hat{\phi}(\vec{z}).
\] 
\end{Lemma1}

\begin{proof}
Let $\theta(x, \vec{z})$ be $\Pi_n$ such that $\phi(\vec{z})$ is $\exists x \theta(x, \vec{z})$. Let $q \in \omega$ be such that $q= \ulcorner \neg \theta(\vec{z}) \urcorner$. Let $z_0, \ldots, z_{m-1} \in M$. Then 
\[
\mathcal{M} \models	 \phi(z_0, \ldots, z_{m-1}) \textrm{ if and only if }  \mathbb{C}\mathrm{ov}^n_\mathcal{M} \models \exists x \exists z(z= \langle x, z_0, \ldots, z_{m-1} \rangle \land \neg \mathsf{S}(q, z)).
\]
\Square
\end{proof}

\begin{Theorems1} \label{Th:MainEndExtensionTheorem}
Let $S$ be a recursively enumerable $\mathcal{L}$-theory such that
\[
S \vdash \mathsf{KP}+\Pi_n\textsf{-Collection}+\Sigma_{n+1}\textsf{-Foundation},
\]
and let $\mathcal{M}= \langle M, \mathsf{E}^\mathcal{M} \rangle$ be a countable model of $S$. Then there exists an $\mathcal{L}$-structure $\mathcal{N}= \langle N, \mathsf{E}^{\mathcal{N}} \rangle$ such that $\mathcal{M} \prec_{e, n} \mathcal{N} \models S$ and there exists $d \in N$ such that for all $x \in M$, $\mathcal{N} \models (x \in d)$. 
\end{Theorems1}

\begin{proof}
Let $T$ be the $\mathcal{L}^{\mathtt{ee}}_{\mathbb{C}\mathrm{ov}^n_\mathcal{M}}$-theory that contains:
\begin{itemize}
\item $S$;
\item for all $a, b \in M$ with $\mathcal{M} \models (a \in b)$, $\bar{a} \in \bar{b}$;
\item for all $a \in M$, 
\[
\forall x \left(x \in a \iff \bigvee_{b \in a} (x=\bar{b}) \right); 
\]
\item for all $a \in M$, $\bar{a} \in \mathbf{c}$;
\item for all $\Pi_n$-formulae, $\phi(x_0, \ldots, x_{m-1})$, and for all $a_0, \ldots, a_{m-1} \in M$ such that $\mathcal{M} \models \phi(a_0, \ldots, a_{m-1})$,
\[
\phi(\bar{a}_0, \ldots, \bar{a}_{m-1}).
\]
\end{itemize}
Since $\mathsf{S}$ is a satisfaction class for $\Sigma_n$-formulae (and hence $\Pi_n$-formula) of $\mathcal{M}$ in $\mathbb{C}\mathrm{ov}^n_\mathcal{M}$, $T \subseteq \mathbb{C}\mathrm{ov}^n_\mathcal{M}$ is $\Sigma_1(\mathcal{L}^*_{\mathsf{S}})$ over $\mathbb{C}\mathrm{ov}^n_\mathcal{M}$. Let $T_0 \subseteq T$ be such that $T_0 \in \mathbb{C}\mathrm{ov}^n_\mathcal{M}$. Using Lemma \ref{Th:KeyCompactnessLemma}, let $c \in M$ be such that
\[
c^*= \{ a \in M \mid \bar{a} \textrm{ is mentioned in } T_0\}.
\]  
Interpreting each $\bar{a}$ that is mentioned in $T_0$ by $a \in M$ and interpreting $\mathbf{c}$ by $c$, we expand $\mathcal{M}$ to a model of $T_0$. Therefore, by the Barwise Compactness Theorem, there exists $\mathcal{N} \models T$. The $\mathcal{L}$-reduct of $\mathcal{N}$ is the desired extension of $\mathcal{M}$. \Square
\end{proof}

\section[Well-founded models of collection]{Well-founded models of collection} \label{sec:application}

In this section we use Theorem \ref{Th:MainEndExtensionTheorem} to show that for all $n \geq 1$, $\mathsf{M}+\Pi_n\textsf{-Collection}+\Pi_{n+1}\textsf{-Foundation}$ proves $\Sigma_{n+1}\textsf{-Separation}$. In particular, the theories $\mathsf{M}+\Pi_n\textsf{-Collection}$ and $\mathsf{M}+\textsf{Strong } \Pi_n\textsf{-Collection}$ have the same well-founded models.

In order to be able to apply Theorem \ref{Th:MainEndExtensionTheorem} to countable models of $\mathsf{M}+\Pi_n\textsf{-Collection}+\Pi_{n+1}\textsf{-Foundation}$, we first need to show that $\mathsf{M}+\Pi_n\textsf{-Collection}+\Pi_{n+1}\textsf{-Foundation}$ proves $\Sigma_{n+1}\textsf{-Foundation}$. The proof presented here generalises the argument presented in \cite[Section 3]{em22} showing that $\mathsf{KP}^{\mathcal{P}}$ proves $\Sigma_1^{\mathcal{P}}\textsf{-Foundation}$.

\begin{Definitions1} \label{Df:FamilyOfPathsThroughPhi}
Let $\phi(x, y, \vec{z})$ be an $\mathcal{L}$-formula. Define $\delta^\phi(a, b, f)$ to be the $\mathcal{L}$-formula:
\[
\begin{array}{c}
(a \in \omega) \land (f \textrm{ is a function}) \land \mathsf{dom}(f)=a+1 \land f(0)= \{b\} \land\\
(\forall u \in \omega) \left(\begin{array}{c}
(\forall x \in f(u))(\exists y \in f(u+1)) \phi(x, y, \vec{z})\\
(\forall y \in f(u+1))(\exists x \in f(u)) \phi(x, y, \vec{z})
\end{array}\right) 
\end{array}
\]
Define $\delta_\omega^\phi(b, f, \vec{z})$ to be the $\mathcal{L}$-formula:
\[
\begin{array}{c}
(f \textrm{ is a function}) \land \mathsf{dom}(f)= \omega \land f(0)= \{b\} \land\\
(\forall u \in \omega) \left(\begin{array}{c}
(\forall x \in f(u))(\exists y \in f(u+1)) \phi(x, y, \vec{z})\\
(\forall y \in f(u+1))(\exists x \in f(u)) \phi(x, y, \vec{z})
\end{array}\right) 
\end{array}
\] 
\end{Definitions1}

Viewing $\vec{z}$ as parameters and letting $a \in \omega$, $\delta^\phi(a, b, f, \vec{z})$ says that $f$ describes a family of directed paths of length $a+1$ starting at $b$ through the directed graph defined by $\phi(x, y, \vec{z})$. Similarly, viewing $\vec{z}$ as parameters, $\delta_\omega^\phi(b, f, \vec{z})$ says that $f$ describes a family of directed paths of length $\omega$ starting at $b$ through the directed graph defined by $\phi(x, y, \vec{z})$. Note that if $\phi(x, y, \vec{z})$ is $\Delta_0$, then, in the theory $\mathsf{M}^-$, both $\delta^\phi(a, b, f \vec{z})$ and $\delta_\omega^\phi(b, f, \vec{z})$ can be written as a $\Delta_0$-formulae with parameter $\omega$. Moreover, if $n \geq 1$ and $\phi(x, y, \vec{z})$ is a $\Sigma_n$-formula ($\Pi_n$-formula), then, in the theory $\mathsf{M}^-+\Pi_{n-1}\textsf{-Collection}$, both $\delta^\phi(a, b, f \vec{z})$ and $\delta_\omega^\phi(b, f, \vec{z})$ can be written as a $\Sigma_n$-formulae ($\Pi_n$-formulae, respectively) with parameter $\omega$.

The following generalises Rathjen's $\Delta_0$-weak dependent choices scheme from \cite{rat92}:
\begin{itemize}
\item[]($\Delta_0\textrm{-}\mathsf{WDC}_\omega$) For all $\Delta_0$-formulae, $\phi(x, y, \vec{z})$, 
\[
\forall \vec{z}(\forall x \exists y \phi(x, y, \vec{z}) \Rightarrow \forall w \exists f \delta_\omega^\phi(w, f, \vec{z}));
\]
\end{itemize}
and for all $n \geq 1$,
\begin{itemize}
\item[]($\Delta_n\textrm{-}\mathsf{WDC}_\omega$) for all $\Pi_n$-formulae, $\phi(x, y, \vec{z})$, and for all $\Sigma_n$-formulae, $\psi(x, y, \vec{z})$,
\[
\forall \vec{z}(\forall x \forall y(\phi(x, y, \vec{z}) \iff \psi(x, y, \vec{z})) \Rightarrow (\forall x \exists y \phi(x, y, \vec{z}) \Rightarrow \forall w \exists f \delta_\omega^\phi(w, f, \vec{z}))).
\]
\end{itemize}

The following is based on the proof of \cite[Proposition 3.2]{rat92}:

\begin{Theorems1} \label{Th:WDCImpliesFoundation}
Let $n \in \omega$ with $n \geq 1$. The theory $\mathsf{KP}+\Pi_{n-1}\textsf{-Collection}+\Sigma_{n}\textsf{-Foundation}+\Delta_{n+1}\textrm{-}\mathsf{WDC}_\omega$ proves $\Sigma_{n+1}\textsf{-Foundation}$. 
\end{Theorems1}

\begin{proof}
Let $T$ be the theory $\mathsf{KP}+\Pi_{n-1}\textsf{-Collection}+\Sigma_{n}\textsf{-Foundation}+\Delta_{n+1}\textrm{-}\mathsf{WDC}_\omega$. Assume, for a contradiction, that $\mathcal{M}= \langle M, \in^\mathcal{M} \rangle$ is such that $\mathcal{M} \models T$ and there is an instance of $\Sigma_{n+1}\textsf{-Foundation}$ that is false in $\mathcal{M}$. Let $\phi(x, y, \vec{z})$ be a $\Pi_n$-formula and let $\vec{a} \in M$ be such that
\[
\{x \mid \mathcal{M} \models \exists y \phi(x, y, \vec{a}) \}
\] 
is nonempty and has no $\in$-minimal element. Let $b, d \in M$ be such that $\mathcal{M} \models \phi(b, d, \vec{a})$. Now,
\[
\mathcal{M} \models \forall x \forall u \exists y \exists v (\phi(x, u, \vec{a}) \Rightarrow (y \in x) \land \phi(y, v, \vec{a})).
\]
Therefore, $\mathcal{M} \models \forall x \exists y \theta(x, y, \vec{a})$ where $\theta(x, y, \vec{a})$ is the formula
\[
x= \langle x_0, x_1 \rangle \land y= \langle y_0, y_1 \rangle \land (\phi(x_0, x_1, \vec{a}) \Rightarrow (y_0 \in x_0) \land \phi(y_0, y_1, \vec{a})).
\]
Now, $\theta(x, y, \vec{a})$ is $\Delta_{n+1}^T$. Work inside $\mathcal{M}$. Using $\Delta_{n+1}\textrm{-}\mathsf{WDC}_\omega$, let $f$ be such that $\delta_\omega^\theta(\langle b, d \rangle, f, \vec{a})$. Now, $\Sigma_n\textsf{-Foundation}$ implies that for all $n \in \omega$,
\begin{itemize}
\item[(i)] $f(n) \neq \emptyset$;
\item[(ii)] for all $x \in f(n)$, $x= \langle x_0, x_1 \rangle$ and $\phi(x_0, x_1, \vec{a})$. 
\end{itemize}
Therefore, for all $n \in \omega$,
\[
\begin{array}{c}
(\forall x \in f(n))(\exists y \in f(n+1))(x= \langle x_0, x_1 \rangle \land y= \langle y_0, y_1 \rangle \land y_0 \in x_0) \land\\
(\forall y \in f(n+1))(\exists x \in f(n))(x= \langle x_0, x_1 \rangle \land y= \langle y_0, y_1 \rangle \land y_0 \in x_0)
\end{array}.
\]
Let $B= \mathsf{TC}(\{b\})$. $\textsf{Set-Foundation}$ implies that for all $n \in \omega$, 
\[
(\forall x \in f(n))(x=\langle x_0, x_1 \rangle \land x_0 \in B).
\]
Let
\[
A= \left\{ x \in B \left| (\exists n \in \omega)(\exists z \in f(n))\left(\exists y \in \bigcup z\right)(z= \langle x, y\rangle) \right. \right\},
\]
which is a set by $\Delta_0\textsf{-Separation}$. Now, let $x \in A$. Therefore, there exists $n \in \omega$ and $z \in f(n)$ such that $z= \langle x, x_0 \rangle$. And, there exists $w \in f(n+1)$ such that $w= \langle y, y_0 \rangle$ and $y \in x$. But $y \in A$. So $A$ is a set with no $\in$-minimal element, which is the desired contradiction. 
\Square
\end{proof}

The following refinement of Definition \ref{Df:FamilyOfPathsThroughPhi} will allow us to show that for $n \geq 1$, $\mathsf{M}+\Pi_n\textsf{-Collection}+\Pi_{n+1}\textsf{-Foundation}$ proves $\Delta_{n+1}\textrm{-}\mathsf{WDC}_\omega$. 

\begin{Definitions1}
Let $\phi(x, y, \vec{z})$ be an $\mathcal{L}$-formula. Define $\eta^\phi(a, b, f, \vec{z})$ to the $\mathcal{L}$-formula:
\[
\begin{array}{c}
\delta^\phi(a, b, f, \vec{z}) \land\\
(\forall u \in a) \exists \alpha \exists X \left(\begin{array}{c}
(\alpha \textrm{ is an ordinal}) \land (X= V_\alpha) \land\\
(\forall x \in f(u+1))(x \in X)\\
(\forall y \in X)(\forall x \in f(u))(\phi(x, y, \vec{z}) \Rightarrow y \in f(u+1)) \land\\
(\forall \beta \in \alpha)(\forall Y \in X)\left( \begin{array}{c}
Y= V_\beta \Rightarrow\\
(\exists x \in f(u))(\forall y \in Y) \neg \phi(x, y, \vec{z})
\end{array}\right)
\end{array}\right)
\end{array}
\]
\end{Definitions1}   

The formula $\eta^\phi(a, b, f, \vec{z})$ says that $f$ is a function with domain $a+1$ and for all $u \in a$, $f(u+1)$ is the set of $y \in V_\alpha$ such that there exists $x \in f(u)$ with $\phi(x, y, \vec{z})$ and $\alpha$ is least such that for all $x \in f(u)$, there exists $y \in V_\alpha$ such that $\phi(x, y, \vec{z})$. In the theory $\mathsf{M}+ \Pi_1\textsf{-Collection}+\Pi_2\textsf{-Foundation}$, `$X=V_\alpha$' can be expressed as both a $\Sigma_2$-formula and a $\Pi_2$-formula. If $n \geq 1$ and, for given parameters $\vec{c}$, $\phi(x, y, \vec{c})$ is equivalent to both a $\Sigma_{n+1}$-formula and a $\Pi_{n+1}$-formula, then, in the theory $\mathsf{M}+ \Pi_n\textsf{-Collection}+\Pi_2\textsf{-Foundation}$, $\eta^\phi(a, b, f, \vec{z})$ is equivalent to a $\Sigma_{n+1}$-formula.

\begin{Theorems1} \label{Th:FoundationImpliesWDC}
Let $n \in \omega$ with $n \geq 1$. The theory $\mathsf{M}+\Pi_n\textsf{-Collection}+\Pi_{n+1}\textsf{-Foundation}$ proves $\Delta_{n+1}\textrm{-}\mathsf{WDC}_\omega$.
\end{Theorems1}

\begin{proof}
Work in the theory $\mathsf{M}+\Pi_n\textsf{-Collection}+\Pi_{n+1}\textsf{-Foundation}$. Let $\phi(x, y, \vec{z})$ be a $\Pi_{n+1}$-formula. Let $\vec{a}$, $b$ be sets and let $\theta(x, y, \vec{z})$ be a $\Sigma_{n+1}$-formula such that
\[
\forall x \forall y(\phi(x, y, \vec{a}) \iff \theta(x, y, \vec{a})). 
\]
We begin by claiming that for all $m \in \omega$, $\exists f \eta^\phi(m, b, f, \vec{a})$. Assume, for a contradiction, that this does not hold. Using $\Pi_{n+1}\textsf{-Foundation}$, let $k \in \omega$ be least such that $\neg \exists f \eta^\phi(k, b, f, \vec{a})$. Since $k \neq 0$, there exists a function $g$ with $\mathsf{dom}(g)=k$ and $\eta^\phi(k-1, b, g, \vec{a})$. Consider the class
\[
A= \{ \alpha \in \mathsf{Ord} \mid \forall X(X= V_\alpha \Rightarrow (\forall x \in g(k-1)) (\exists y \in X) \phi(x, y, \vec{a}))\}.
\]
\[
= \{ \alpha \in \mathsf{Ord} \mid \exists X(X=V_\alpha \land (\forall x \in g(k-1))(\exists y \in X) \theta(x, y, \vec{a})) \}
\]
Applying $\Sigma_{n+1}\textsf{-Collection}$ to the formula $\theta(x, y, \vec{a})$ shows that $A$ is nonempty. Moreover, $\Delta_{n+1}\textsf{-Foundation}$ ensures that there is a least element $\beta \in A$. Now, let
\[
C= \{y \in V_\beta \mid (\exists x \in g(k-1)) \phi(x, y, \vec{a})\}, 
\]
which is a set by $\Delta_{n+1}\textsf{-Separation}$. Let $f= g \cup \{\langle k, C\}$. Then $f$ is such that $\eta^\phi(k, b, f, \vec{a})$, which contradicts our assumption that no such $f$ exists. Therefore, for all $m \in \omega$, $\exists f \eta^\phi(m, b, f, \vec{a})$. Using $\Sigma_{n+1}\textsf{-Collection}$, let $D$ be such that $(\forall m \in \omega)(\exists f \in D) \eta^\phi(m, b, f, \vec{a})$. Note that for all $m \in \omega$ and for all functions $f$ and $g$, if $\eta^\phi(m, b, f, \vec{a})$ and $\eta^\phi(m, b, g, \vec{a})$, then $f=g$. Now, let
\[
h= \{\langle m, X \rangle \in \omega \times \mathsf{TC}(D) \mid (\exists f \in D)(\eta^\phi(m, b, f, \vec{a}) \land f(m)= X)\}.
\]
Since
\[
h= \{\langle m, X \rangle \in \omega \times \mathsf{TC}(D) \mid (\forall f \in D)(\eta^\phi(m, b, f, \vec{a}) \Rightarrow f(m)=X) \},
\]
$h$ is a set by $\Delta_{n+1}\textsf{-Separation}$. Now, $h$ is the function required by $\Delta_{n+1}\textsf{-}\mathsf{WDC}_\omega$.      
\Square
\end{proof}

Note $\Pi_{n+1}\textsf{-Foundation}$ is only used in the proof of Theorem \ref{Th:FoundationImpliesWDC} to find the least element of a $\Pi_{n+1}$-definable subclass of naturals numbers. Therefore, the proof of Theorem \ref{Th:FoundationImpliesWDC} also yields the following result. 

\begin{Theorems1} \label{Th:StandardImpliesWDC}
Let $n \in \omega$ with $n \geq 1$. Let $\mathcal{M}$ be an $\omega$-standard model of $\mathsf{M}+\Pi_n\textsf{-Collection}+\Pi_2\textsf{-Foundation}$. Then 
\[
\mathcal{M} \models \Delta_{n+1}\textsf{-}\mathsf{WDC}_\omega.
\]
\Square
\end{Theorems1}

Note that $\Pi_2\textsf{-Foundation}$ coupled with $\Pi_1\textsf{-Collection}$ ensures that the function $\alpha \mapsto V_\alpha$ is total.

Combining Theorem \ref{Th:WDCImpliesFoundation} with Theorems \ref{Th:FoundationImpliesWDC} and \ref{Th:StandardImpliesWDC} yields:

\begin{Coroll1} \label{Th:PiFoundationImpliesSigmaFoundation}
Let $n \in \omega$ with $n \geq 1$. The theory $\mathsf{M}+\Pi_n\textsf{-Collection}+\Pi_{n+1}\textsf{-Foundation}$ proves $\Sigma_{n+1}\textsf{-Foundation}$. \Square
\end{Coroll1}

\begin{Coroll1} \label{Th:StandardImpliesSigmaFoundation}
Let $n \in \omega$ with $n \geq 2$. Let $\mathcal{M}$ be an $\omega$-standard model of $\mathsf{M}+\Pi_n\textsf{-Collection}$. Then 
\[
\mathcal{M} \models \Sigma_{n+1}\textsf{-Foundation}.
\]
\Square  
\end{Coroll1}

The proof of \cite[Theorem 3.11]{em22} shows how the use of the cumulative hierarchy can be avoided in the argument used in the proof of Theorem \ref{Th:FoundationImpliesWDC}. The following is \cite[Corollary 3.12]{em22} combined with \cite[Proposition Scheme 6.12]{mat01} and provides a version of Corollary \ref{Th:StandardImpliesSigmaFoundation} when $n=1$:

\begin{Theorems1} \label{Th:BaseCaseStandardImpliesSigmaFoundation}
Let $\mathcal{M}$ be an $\omega$-standard model of $\mathsf{MOST}+\Pi_1\textsf{-Collection}$. Then 
\[
\mathcal{M} \models \Sigma_{2}\textsf{-Foundation}.
\]
\Square  
\end{Theorems1}

Equipped with these results, we are now able to show that, in the theory $\mathsf{M}+\Pi_n\textsf{-Collection}$, $\Pi_{n+1}\textsf{-Foundation}$ implies $\Sigma_{n+1}\textsf{-Separation}$.

\begin{Lemma1}
Let $\mathcal{M}= \langle M, \in^\mathcal{M} \rangle$ and $\mathcal{N}=\langle N, \in^\mathcal{N} \rangle$ be such that $\mathcal{M}, \mathcal{N} \models \mathsf{M}$. If $\mathcal{M} \prec_{e, 1} \mathcal{N}$, then $\mathcal{M} \subseteq_e^\mathcal{P} \mathcal{N}$.
\end{Lemma1}

\begin{proof}
Assume that $\mathcal{M} \prec_{e, 1} \mathcal{N}$. Let $x \in M$ and let $y \in N$ with $\mathcal{N} \models (y \subseteq x)$. We need to show that $y \in M$. Let $a \in M$ be such that $\mathcal{M} \models (a= \mathcal{P}(x))$. Therefore, $\mathcal{M} \models \theta(x, a)$ where $\theta(x, a)$ is the $\Pi_1$-formula
\[
\forall z(z \subseteq x \iff z \in a).
\]
So, $\mathcal{N} \models \phi(x, a)$. Therefore, $\mathcal{N} \models (y \in a)$ and so $y \in N$. \Square
\end{proof} 

As alluded to in \cite[Remark 3.21]{mat01}, the theory $\mathsf{KP}+\Sigma_1\textsf{-Separation}$ is capable of endowing any well-founded partial order with a ranking function.

\begin{Lemma1} \label{Th:WFPartialOrdersHaveRanks}
The theory $\mathsf{KP}+\Sigma_1\textsf{-Separation}$ proves that if $\langle X, R \rangle$ is a well-founded strict partial order, then there exists an ordinal $\gamma$ and a function $h: X \longrightarrow \gamma$ such that for all $x, y \in X$, if $\langle x, y \rangle \in R$, then $h(x) < h(y)$.
\end{Lemma1}

\begin{proof}
Work in the theory $\mathsf{KP}+\Sigma_1\textsf{-Separation}$. Let $X$ be a set and $R \subseteq X \times X$ be such that $\langle X, R \rangle$ is a well-founded strict partial order. Let $\theta(x, g, X, R)$ be the conjunction of the following clauses:
\begin{itemize}
\item[(i)] $g$ is a function;
\item[(ii)] $\mathsf{rng}(g)$ is a set of ordinals;
\item[(iii)] $\mathsf{dom}(g)= \{y \in X \mid \langle y, x \rangle \in R \lor y= x\}$;
\item[(iv)] $(\forall y, z \in \mathsf{dom}(g))(\langle y, z \rangle \in R \Rightarrow g(y) < g(z))$;
\item[(v)] $(\forall y \in \mathsf{dom}(g))(\forall \alpha \in g(y))(\exists z \in X)(\langle z, y \rangle \in R \land g(z) \geq \alpha)$. 
\end{itemize}
Note that $\theta(x, g, X, R)$ can be written as a $\Delta_0$-formula. Moreover, for all $x \in X$ and functions $g_0$ and $g_1$, if $\theta(x, g_0, X, R)$ and $\theta(x, g_1, X, R)$, then $g_0=g_1$. And, if $x, y \in X$ with $\langle x, y \rangle \in R$ and $g_0$ and $g_1$ are functions with $\theta(y, g_0, X, R)$ and $\theta(x, g_1, X, R)$, then $g_0= g_1 \upharpoonright \mathsf{dom}(g_0)$. Now, consider 
\[
A= \{ x \in X \mid \neg \exists g \theta(x, g, X, R)\}, 
\] 
which is a set by $\Pi_1\textsf{-Separation}$. Assume, for a contradiction, that $A \neq \emptyset$. Let $x_0 \in A$ be $R$-minimal. Let $B= \{y \in X \mid \langle y, x_0 \rangle \in R\}$. Using $\Delta_0\textsf{-Collection}$, let $C_0$ be such that $(\forall y \in B)(\exists g \in C_0) \theta(y, g, X, R)$. Let
\[
D= \{ g \in C_0 \mid (\exists y \in B) \theta(y, g, X, R)\}.
\]
Let 
\[
\beta= \sup\{g(y)+1 \mid y \in B \textrm{ and } g \in D \textrm{ with } y \in \mathsf{dom}(g)\}.
\]  
Then $f= \bigcup D \cup \{ \langle x_0, \beta \rangle\}$ is such that $\theta(x_0, f, X, R)$, which contradicts the fact that $x_0 \in A$. Therefore, $A= \emptyset$. Using $\Delta_0\textsf{-Collection}$, let $C_1$ be such that $(\forall x \in X)(\exists g \in C_1) \theta(x, g, X, R)$. Let 
\[
F= \{ g \in C_1 \mid (\exists x \in X) \theta(x, g, X, R)\}.
\]
Then $h= \bigcup F$ is the function we require.      
\Square
\end{proof} 

\begin{Theorems1} \label{Th:FoundationImpliesSeparation}
Let $n \in \omega$ with $n \geq 1$. The theory $\mathsf{M}+\Pi_n\textsf{-Collection}+\Pi_{n+1}\textsf{-Foundation}$ proves $\Sigma_{n+1}\textsf{-Separation}$. 
\end{Theorems1} 

\begin{proof}
Let $\mathcal{M}= \langle M, \in^\mathcal{M} \rangle$ be such that $\mathcal{M} \models \mathsf{M}+\Pi_n\textsf{-Collection}+\Pi_{n+1}\textsf{-Foundation}$. Let $\theta(x, y, \vec{z})$ be a $\Pi_n$-formula and let $b, \vec{a} \in M$. We need to show that $A= \{x \in b \mid \exists y \theta(x, y, \vec{a})\}$ is a set in $\mathcal{M}$. By Corollary \ref{Th:PiFoundationImpliesSigmaFoundation}, $\mathcal{M} \models \Sigma_{n+1}\textsf{-Foundation}$. Using Theorem \ref{Th:MainEndExtensionTheorem}, let $\mathcal{N}= \langle N, \in^\mathcal{N} \rangle$ be such that $\mathcal{M} \prec_{e, n} \mathcal{N}$, $\mathcal{N} \models \mathsf{M}+\Pi_n\textsf{-Collection}+\Pi_{n+1}\textsf{-Foundation}$ and there exists $d \in N$ such that for all $x \in M$, $\mathcal{N} \models (x \in d)$. Let $\alpha \in \mathsf{Ord}^\mathcal{N}$ be such that for all $x \in M$, $\mathcal{M} \models (x \in V_\alpha)$. 

Work inside $\mathcal{N}$. Let
\[
D= \{x \in b \mid (\exists y \in V_\alpha) \theta(x, y, \vec{a})\},
\]
which is a set by $\Pi_n\textsf{-Separation}$. Let
\[
g= \left\{\langle x, \beta \rangle \in D \times \alpha \left| \begin{array}{c}
(\exists y \in V_\alpha)(\rho(y)= \beta \land \theta(x, y, \vec{a})) \land\\
(\forall z \in V_\alpha)(\phi(x, z, \vec{a}) \Rightarrow \beta \leq \rho(z))
\end{array} \right\}\right.,
\]
which is a set by $\Delta_{n+1}\textsf{-Separation}$. Moreover, $g$ is a function. Let $\lhd= \{ \langle x_0, x_1 \rangle \in D \times D \mid g(x_0) < g(x_1)\}$. Note that $\lhd$ is a well-founded strict partial order on $D$.

Since $\mathcal{M} \subseteq_e^\mathcal{P} \mathcal{N}$, $D, \lhd \in M$. Moreover,
\[
\mathcal{M} \models (\lhd \textrm{ is a well-founded strict partial order on } D).
\]

Work inside $\mathcal{M}$. Since $\mathcal{M} \prec_{e, n} \mathcal{N}$, for all $x \in b$, if $\exists y \theta(x, y, \vec{a})$, then $x \in D$. And, for all $x_0, x_1 \in D$, if $\exists y \theta(x_0, y, \vec{a})$ and $\neg \exists y \theta(x_1, y, \vec{a})$, then $x_0 \lhd x_1$. Using Lemma \ref{Th:WFPartialOrdersHaveRanks}, let $\gamma$ be an ordinal and let $h: D \longrightarrow \gamma$ be such that for all $x_0, x_1 \in D$, if $\langle x_0, x_1 \rangle \in D$, then $h(x_0) < h(x_1)$. Consider the class
\[
B= \{ \beta \in \gamma \mid (\exists x \in D)(h(x)= \beta \land \neg \exists y \theta(x, y, \vec{a}))\}.
\]
If $B$ is empty, then $D= \{x \in b \mid \exists y \phi(x, y, \vec{a})\}$ and we are done. Therefore, assume that $B$ is nonempty. So, by $\Pi_{n+1}\textsf{-Foundation}$, $B$ has a least element $\xi$. Let $D_\xi= \{x \in D \mid h(x) < \xi \}$. Let $x \in D_\xi$. Since $\xi$ is the least element of $B$ and $h(x) < \xi$, $\exists y \theta(x, y, \vec{a})$. Conversely, let $x \in b$ be such that $\exists y \theta(x, y, \vec{a})$. Let $x_0 \in D$ be such $h(x_0)= \xi$ and $\neg \exists y \theta(x_0, y, \vec{a})$. Since $\exists y \theta(x, y, \vec{a})$, it must be the case that $h(x) < h(x_0)=\xi$. So, $x \in D_\xi$. This shows that $D_\xi=\{x \in b \mid \exists y \theta(x, y, \vec{a})\}$. Therefore, $\Sigma_{n+1}\textsf{-Separation}$ holds in $\mathcal{M}$. \Square     
\end{proof}

Gostanian \cite{gos80} notes that the techniques he uses to compare the heights of minimum models of subsystems of $\mathsf{ZF}$ without the powerset axiom do not apply to subsystems that include the powerset axiom. Theorem \ref{Th:FoundationImpliesSeparation} settles the relationship between the heights of the minimum models of the theories $\mathsf{M}+\Pi_n\textsf{-Collection}$ and $\mathsf{M}+\textsf{Strong } \Pi_n\textsf{-Collection}$ for all $n \geq 1$.

\begin{Coroll1}
Let $n \in \omega$ with $n \geq 1$. The theories $\mathsf{M}+\Pi_n\textsf{-Collection}$ and $\mathsf{M}+\textsf{Strong } \Pi_n\textsf{-Collection}$ have the same transitive models. In particular, the minimum models $\mathsf{M}+\Pi_n\textsf{-Collection}$ and $\mathsf{M}+\textsf{Strong } \Pi_n\textsf{-Collection}$ coincide. \Square 
\end{Coroll1}

The results of \cite{mck19} show that for all $n \geq 1$, $\mathsf{M}+\textsf{Strong }\Pi_n\textsf{-Collection}$ proves the consistency of $\mathsf{M}+\Pi_n\textsf{-Collection}$. Theorem \ref{Th:FoundationImpliesSeparation} yields the following:

\begin{Coroll1}
Let $n \in \omega$ with $n \geq 1$. The theory $\mathsf{M}+\textsf{Strong } \Pi_n\textsf{-Collection}$ does not prove the existence of a transitive model of $\mathsf{M}+\Pi_n\textsf{-Collection}$. \Square
\end{Coroll1}  

The following example shows that the statement of Theorem \ref{Th:FoundationImpliesSeparation} with $n=0$ does not hold.

\begin{Examp1}
Let $\mathcal{M}= \langle M, \in^{\mathcal{M}} \rangle$ be an $\omega$-standard model of $\mathsf{ZFC}$ in which there is a countable ordinal that is nonstandard. Note that such a model can built from a transitive model of $\mathsf{ZFC}$ using, for example, \cite[Theorem 2.4]{km68}, or using the Barwise Compactness Theorem as in \cite[Lemma 7.2]{mck15}. Let $W$ be the transitive set that is isomorphic to the well-founded part of $\mathcal{M}$. Then $\langle W, \in \rangle$ satisfies $\mathsf{KP}^\mathcal{P}+\textsf{Foundation}$. However, there are well-orderings of $\omega$ in $\langle W, \in \rangle$ that are not isormorphic to any ordinal in $\langle W, \in \rangle$, so $\langle W, \in \rangle$ does not satisfy $\Sigma_1\textsf{-Separation}$. \Square  
\end{Examp1}

The following is a consequence of Theorems 2.1 and 2.2 in \cite{gos80} and shows that the presence of $\textsf{Powerset}$ is essential in Theorem \ref{Th:FoundationImpliesSeparation}.

\begin{Theorems1}
(Gostanian) Let $n \in \omega$. Let $\alpha$ be the least ordinal such that $\langle L_\alpha, \in \rangle \models \mathsf{KP}+\Pi_n\textsf{-Collection}$. Then $\langle L_\alpha, \in \rangle$ does not satisfy $\Sigma_{n+1}\textsf{-Separation}$. \Square 
\end{Theorems1}    

In \cite[Theorem 4.6]{res87} (see also \cite[Theorem 4.15]{flw16}), Ressayre shows that for all $n \in \omega$, the theory $\mathsf{KP}+\mathsf{V}=\mathsf{L}+\Pi_n\textrm{-Collection}+\Sigma_{n+1}\textsf{-Foundation}$ does not prove $\Pi_{n+1}\textsf{-Foundation}$. Ressayre's construction can be adapted (as noted in \cite[Theorem 4.15]{res87}) to show that for all $n \geq 1$, $\mathsf{M}+\Pi_n\textsf{-Collection}+\Sigma_{n+1}\textsf{-Foundation}$ does not prove $\Pi_{n+1}\textsf{-Foundation}$. Since $\mathsf{M}+\Sigma_{n+1}\textsf{-Separation}$ proves, $\Pi_{n+1}\textsf{-Foundation}$, this shows that $\mathsf{M}+\Pi_n\textsf{-Collection}+\Sigma_{n+1}\textsf{-Foundation}$ does not prove $\Sigma_{n+1}\textsf{-Separation}$.

\begin{Theorems1}
(Ressayre) Let $n \in \omega$ with $n \geq 1$. The theory $\mathsf{M}+\Pi_n\textsf{-Collection}+\Sigma_{n+1}\textsf{-Foundation}$ does not prove $\Pi_{n+1}\textsf{-Foundation}$.
\end{Theorems1}

\begin{proof}
Let $\mathcal{M}= \langle M, \in^\mathcal{M} \rangle$ be a nonstandard $\omega$-standard model of $\mathsf{ZF}+\mathsf{V}=\mathsf{L}$. Let $\delta \in \mathsf{Ord}^\mathcal{M}$ be nonstandard. Let $I \subseteq (\delta+\delta)^*$ be an initial segment of $(\delta+\delta)^*$ such that $(\delta+\delta)^* \backslash I$ has no least element.

Work inside $\mathcal{M}$. Define a function $f$ with domain $\delta+\delta$ such that
\[
\begin{array}{ll}
f(0)= V_\gamma & \textrm{where } \gamma \textrm{ is least such that } V_\gamma \textrm{ is a } \Sigma_n\textrm{-elementary}\\
&\textrm{substructure of the universe};\\
f(\alpha+1)= V_\gamma & \textrm{where } \gamma \textrm{ is least such that } f(\alpha) \in V_\gamma \textrm{ and}\\
&V_\gamma \textrm{ is a } \Sigma_n\textrm{-elementary substructure of the universe};\\
f(\beta)= \bigcup_{\alpha \in \beta} f(\alpha) & \textrm{if } \beta \textrm{ is a limit ordinal}.
\end{array}
\]

Now, working in the metatheory again, define $\mathcal{N}= \langle N, \in^\mathcal{N} \rangle$ by:
\[
N= \bigcup_{\alpha \in I} f(\alpha)^* \textrm{ and } \in^\mathcal{N} \textrm{ is the restriction of } \in^\mathcal{M} \textrm{ to } N.
\]
Therefore, $\mathcal{N} \prec_{e, n} \mathcal{M}$ and $\mathsf{Ord}^\mathcal{M} \backslash \mathsf{Ord}^\mathcal{N}$ has no least element. It is clear that $\mathcal{N}$ is $\omega$-standard and satisfies $\mathsf{M}+\mathsf{AC}$. We claim that $\mathcal{N}$ satisfies $\textsf{Strong } \Delta_0\textsf{-Collection}$. Let $\phi(x, y, \vec{z})$ be a $\Delta_0$-formula, and let $b, \vec{a} \in N$. Let $\alpha \in \mathsf{Ord}^\mathcal{N}$ be such that $V_\alpha^\mathcal{M} \in N$, $b, \vec{a} \in (V_\alpha^\mathcal{M})^*$ and $\langle (V_\alpha^\mathcal{M})^*, \in^\mathcal{N} \rangle \prec_{e, 1} \mathcal{N}$. But then
\[
\mathcal{N} \models (\forall x \in b)(\exists y \phi(x, y, \vec{a}) \Rightarrow (\exists y \in V_\alpha) \phi(x, y, \vec{a})).
\]
This shows that $\mathcal{N}$ satisfies $\textsf{Strong } \Delta_0\textsf{-Collection}$. So, $\mathcal{N} \models \mathsf{MOST}+\mathsf{V}=\mathsf{L}$. Therefore, by Theorem \ref{Th:CollectionFromPartiallyElementaryEE1},
\[
\mathcal{N} \models \mathsf{MOST}+ \Pi_n\textsf{-Collection}.
\]
And, by Theorem \ref{Th:BaseCaseStandardImpliesSigmaFoundation} ($n=1$) and Corollary \ref{Th:StandardImpliesSigmaFoundation} ($n > 1$), 
\[
\mathcal{N} \models \Sigma_{n+1}\textsf{-Foundation}.
\]
Note that `$X$ is $\Sigma_n$-elementary submodel of the universe', which we abbreviate $X \prec_n \mathbb{V}$, can be expressed as
\[
(\forall x \in X^{<\omega})(\forall m \in \omega)(\mathsf{Sat}_{\Sigma_n}(m, x) \Rightarrow \langle X, \in \rangle \models \mathsf{Sat}_{\Sigma_n}(m, x)),
\]
and is equivalent to a $\Pi_n$-formula. Now, consider the formula $\theta(\alpha)$ defined by
\[
\exists f \left(\begin{array}{c}
(f \textrm{ is a function}) \land \mathsf{dom}(f)= \alpha \land\\
\exists X \exists \beta (
X= V_\beta \land X \prec_n \mathbb{V} \land f(0)=X \land (\forall Y, \gamma \in X)(Y= V_\gamma \Rightarrow \neg (Y \prec_n \mathbb{V})))\land\\
(\forall \eta \in \alpha)\left( \begin{array}{c}
\eta= \xi+1 \Rightarrow \exists X \exists \beta \left(\begin{array}{c}
X=V_\beta \land X \prec_n \mathbb{V} \land f(\eta)= X \land f(\xi) \in X \land\\
(\forall Y, \gamma \in X)(Y \neq V_\gamma \lor \neg(Y \prec_n \mathbb{V}) \lor f(\xi) \notin Y)
\end{array}\right)
\end{array}\right)\\
\land (\forall \eta \in \alpha)\left( \begin{array}{c}
(\eta \textrm{ is a limit}) \Rightarrow f(\eta)= \bigcup_{\xi \in \eta} f(\xi)
\end{array} \right)
\end{array} \right).
\]
Note that $\theta(\alpha)$ can be expressed as a $\Sigma_{n+1}$-formula and says that there exists a function that enumerates the first $\alpha$ levels of the cumulative hierarchy that are $\Sigma_n$-elementary submodels of the universe. Therefore, the class
\[
A= \{ \alpha \in \mathsf{Ord}^\mathcal{N} \mid \neg \theta(\alpha) \}= \mathsf{Ord}^\mathcal{N} \backslash I
\]
has no least element, so $\Pi_{n+1}\textsf{-Foundation}$ fails in $\mathcal{N}$.
\Square  
\end{proof}

\bibliographystyle{alpha}
\bibliography{.}                 

\begin{thebibliography}{9}

\bibitem[Bar74]{bar74} Barwise, K. J. ``Admissible sets over models of set theory". {\it Generalized Recursion Theory: Proceedings of the 1972 Oslo Symposium}. Edited by J.F. Fenstad and P. G. Hinman. North-Holland, Amsterdam. 1974.

\bibitem[Bar75]{bar75} Barwise, K. J. {\it Admissible Sets and Structures}. Perspectives in Mathematical Logic. Springer-Verlag, Berlin-Heidelberg-New York. 1975.

\bibitem[EM22]{em22} Enayat, A. and McKenzie, Z. {\it End extending models of set theory via power admissible covers}. Annals of Pure and Applied Logic. Vol. 173. No. 8. 2022. pp 103--132.

\bibitem[Fri]{fri73} Friedman, H. M. ``Countable models of set theories". \emph{Cambridge Summer School in Mathematical Logic, August 1--21, 1971}. Edited by A. R. D. Mathias and H. Rogers Jr. Springer Lecture Notes in Mathematics. Vol. 337. Springer, Berlin. 1973. pp 539--573.

\bibitem[FLW]{flw16} Friedman, S.-D., Li, W.; and Wong, T. L. ``Fragments of Kripke-Platek Set Theory and the Metamathematics of $\alpha$-Recursion Theory". \emph{Archive for Mathematical Logic}. Vol. 55. No. 7. 2016. pp 899--924.

%\bibitem[Hut]{hut76} Hutchinson, John E. ``Elementary Extensions of Countable Models of Set Theory". \emph{The Journal of Symbolic Logic}. Vol. 41. No. 1. 1976. pp 139--145.

\bibitem[Gos]{gos80} Gostanian, R. ``Constructible models of subsystems of $\mathrm{ZF}$". {\it The Journal of Symbolic Logic}. Vol. 45. No. 2. 1980. pp 237--250.

%\bibitem[JS]{js18} J\"{a}ger, G.; and Steila, S. ``About some fixed point axioms and related principles in Kripke-Platek environments". {\it The Journal of Symbolic Logic}. Vol. 83. No. 2. 2018. pp 642--668.

\bibitem[Kau]{kau81} Kaufmann, M. ``On Existence of $\Sigma_n$ End Extensions". \emph{Logic Year 1979-80, The University of Connecticut}. Lecture Notes in Mathmeatics. No. 859. Springer-Verlag. 1981. pp 92--103.

\bibitem[KM]{km68} Keisler, H. J.; and Morley, M. ``Elementary extensions of models of set theory". {\it Israel Journal of Mathematics}. Vol. 5. 1968. pp 49--65

%\bibitem[MS]{ms61} MacDowell, R.; and Specker, E. ``Modelle der Arithmetik" in {\it Infinitistic Methods}. Proceedings of the Symposium on the Foundations of Mathematics, Warsaw, September 1959. Pergamon Press, Oxford and Panstowa Wydanictwo Naukowe, Warsaw. 1961. pp 257--263.

\bibitem[M15]{mck15} McKenzie, Z. ``Automorphisms of models of set theory and extensions of NFU". {\it Annals of Pure and Applied Logic}. Vol. 166. 2015. pp 601--638. 

\bibitem[M19]{mck19} McKenzie, Z. ``On the relative strengths of fragments of collection". {\it Mathematical Logic Quarterly}. Vol. 65. No. 1. 2019. pp 80--94.

%\bibitem[Mat69]{mat69} Mathias, Adrian R. D. ``Notes on set theory". Available online: {\tt http://www.dpmms.cam.ac.uk/\textasciitilde ardm/} (last accessed on 29/iv/2018)

\bibitem[Mat01]{mat01} Mathias, A. R. D. ``The strength of Mac Lane set theory". \emph{Annals of Pure and Applied Logic}. Vol. 110. 2001. pp 107-234.

%\bibitem[PK]{pk78} Paris, J. B.; and Kirby, L. A. S. ``$\Sigma_n$-collection schemas in arithmetic". In {\it Logic Colloquium '77 (Proceedings of the colloquium held in Wroc\l aw, August 1977)}. {\it Studies in Logic and the Foundations of Mathematics}. Vol. 96. North-Holland, Amsterdam-New York, 1978. pp 199--209.

\bibitem[Rat92]{rat92} Rathjen, Michael. ``A proof-theoretic characterization of the primitive recursive set functions". {\it The Journal of Symbolic Logic}. Vol. 57. No. 3. 1992. pp 954--969.

\bibitem[Res]{res87} Ressayre, J.-P. ``Mod\`{e}les non standard et sous-syst\`{e}mes remarquables de ZF". In {\it Mod\`{e}les non standard en arithm\'{e}tique et th\'{e}orie des ensembles}. Volume 22 of {\it Publications Math\'{e}matiques de l'Universit\'{e} Paris VII}. Universit\'{e} de Paris VII, U.E.R. de Math\'{e}matiques, Paris, 1987. pp 47--147. 

\bibitem[Tak]{tak72} Takahashi, M. ``$\tilde{\Delta}_1$-definability in set theory". \emph{Conference in mathematical logic --- London '70}. Edited by W. Hodges. Springer Lecture Notes in Mathematics. Vol. 255. Springer. 1972. pp 281-304.

\end{thebibliography}

\end{document}